\documentstyle[11pt,amssymb,pb-diagram]{article}

\newlength{\minitwocolumn}
\newcommand{\Z}{{\Bbb Z}} 
\newcommand{\R}{{\Bbb R}} 
\newcommand{\C}{{\Bbb C}} 
\newcommand{\F}{{\cal F}}

\newcommand{\cR}{{\cal R}}

\newcommand{\hL}{\widehat{L}}

\renewcommand{\H}{{\cal H}}
\newcommand{\la}{\lambda}

\newcommand{\nn}{{\nonumber}}
\newcommand{\eqref}[1]{(\ref{#1})}
\newcommand{\bea}{\begin{eqnarray}}
\newcommand{\ena}{\end{eqnarray}}
\newcommand{\beit}{\begin{itemize}}
\newcommand{\enit}{\end{itemize}}
\newcommand{\be}{\begin{eqnarray*}}
\newcommand{\en}{\end{eqnarray*}}
\newcommand{\lb}[1]{\label{#1}}

\newcommand{\id}{\hbox{id}}

\newcommand{\lal}{\langle}
\newcommand{\ral}{\rangle}

\def\infq4p#1{{(#1;q^4,p)_\infty}}

\newcommand{\al}{\alpha}

\newcommand{\ep}{\epsilon}
\newcommand{\vep}{\varepsilon}
\newcommand{\bep}{\bar{\epsilon}}
\newcommand{\bbep}{\eta}
\newcommand{\ba}{\bar{\alpha}}

\newcommand{\bfv}{{\bf v}}
\newcommand{\du}[1]{{\frac{d {#1}}{2\pi i {#1}}}}

\font\teneufm=eufm10
\font\seveneufm=eufm7
\font\fiveeufm=eufm5
\newfam\eufmfam
\textfont\eufmfam=\teneufm
\scriptfont\eufmfam=\seveneufm
\scriptscriptfont\eufmfam=\fiveeufm

\let\goth\frak
\newcommand{\slth}{\widehat{\goth{sl}}_2}
\newcommand{\slt}{\goth{sl}_2}
\newcommand{\slnh}{\widehat{\goth{sl}}_N}
\newcommand{\sln}{\goth{sl}_N}
\newcommand{\g}{\goth{g}}

\newcommand{\Aqp}{{\cal A}_{q,p}}

\newcommand{\Bqla}{{{\cal B}_{q,\lambda}}}

\newcommand{\slnhbig}{\widehat{\mbox{\fourteeneufm sl}}_N}  

\newcommand{\h}{\goth{h}}

\font\fourteeneufm=eufm10 scaled\magstep2    
   
\makeatletter
\@addtoreset{equation}{section}
\makeatother

\newtheorem{thm}{Theorem}[section]
\newtheorem{prop}[thm]{Proposition}
\newtheorem{lem}[thm]{Lemma}
\newtheorem{cor}[thm]{Corollary}
\newtheorem{conj}[thm]{Conjecture}

\newtheorem{df}{Definition}[section]
\newtheorem{dfn}[thm]{Definition}

\begin{document}

\

\vspace{0.7cm}
\begin{center}
{\Large{\bf
The Elliptic Algebra $U_{q,p}(\slnhbig)$ and }}

 
{\Large{\bf the Drinfeld Realization of 
the Elliptic Quantum Group $\Bqla(\slnhbig)$
}}

\vspace{1.2cm}
{ Takeo KOJIMA$~^{*,\dagger}$ and~ Hitoshi KONNO$~^{**,\dagger}$}

\vspace{0.7cm}
{\it
~*
Department of Mathematics,
College of Science and Technology,\\
Nihon University, Chiyoda-ku, Tokyo
101-0062, Japan.\\
E-mail:kojima@math.cst.nihon-u.ac.jp\\~\\
~**
Department of Mathematics,
Faculty of Integrated Arts and Sciences,\\
Hiroshima University,
Higashi-Hiroshima 739-8521, Japan.
\\
E-mail:konno@mis.hiroshima-u.ac.jp
}

\vspace{0.8cm}
{\it~$\dagger$ 
Department of Mathematics,  Heriot-Watt University,\\
 Edinburgh EH14 4AS, UK. }
 
\end{center}

\vspace{0.8cm}
\begin{abstract}
By using the elliptic analogue of the Drinfeld
currents in the elliptic algebra $U_{q,p}(\slnh)$,
 we construct a $L$-operator, which satisfies the 
$RLL$-relations characterizing the face type 
elliptic quantum group $\Bqla(\slnh)$. For this 
purpose, we introduce a set of new currents $K_j(v)\ 
 (1\leq j\leq N)$ in $U_{q,p}(\slnh)$. As in the $N=2$ case,
 we find a structure of $U_{q,p}(\slnh)$ as a certain
 tensor product of $\Bqla(\slnh)$ and a Heisenberg algebra.  
In the level-one representation, we give a
free field realization of the currents in $U_{q,p}(\slnh)$. 
Using the coalgebra structure of $\Bqla(\slnh)$ and 
the above tensor 
structure, we derive a free field
realization of the 
$U_{q,p}(\slnh)$-analogue of 
$\Bqla(\slnh)$-intertwining operators. 
The resultant operators 
coincide with those of the vertex operators in the $A_{N-1}^{(1)}$-type
face model. 
\end{abstract}

\newpage
\section{Introduction}
In recent papers\cite{FIJKMY1,Felder,Fronsdal,EF,JKOS1}, 
the notion of elliptic quantum groups
has been proposed. There are  two types of elliptic quantum groups, the 
vertex type $\Aqp(\slnh)$ 
and the face type $\Bqla(\g)$, where $\g$ is a 
Kac-Moody algebra associated with 
a symmetrizable generalized Cartan matrix.
The elliptic quantum groups have the structure of 
quasi-triangular quasi-Hopf algebras introduced by 
Drinfeld \cite{Drinfeld0}.  
Since certain finite dimensional 
representations of the universal $R$-matrices 
of these elliptic quantum groups
yield known elliptic Boltzmann weights including, for example, those of 
the eight vertex model\cite{Baxter} and  the Andrews-Baxter-Forrester (ABF) face model\cite{ABF}, 
we expect that we can perform an algebraic analysis of both types of 
elliptic lattice models based on 
the corresponding elliptic quantum groups. 

Here, algebraic analysis means, in a restricted sense, 
a method of studying two 
dimensional solvable lattice models  based on the representation theory of 
infinite dimensional quantum groups \cite{JM}. It can be regarded as an
off-critical extension of 
conformal field theory, where the representation theory of    
the Virasoro algebras and/or affine Lie algebras plays an
essential role.
In fact, quite a lot of, but not all, solvable lattice models allow us, 
in the thermodynamic limit, to identify the space of states of the 
models with the infinite dimensional modules of certain quantum groups. 
Then two types of intertwining operators,  type I and
 type II, of such modules become important. The type I intertwiner 
provides a realization of local operators, such as spin operators for example, 
on the infinite dimensional 
modules of quantum groups.  
And the type II plays the  role of creation operator of physical excitations. 
Due to the coalgebra structure of quantum groups, these intertwiners can be  
determined uniquely. Realizing these ingredients in certain forms, such as 
the free field realization for example, 
one can perform a calculation of correlation functions as well 
as form factors of the models.
  
Through  experience of the analysis of 
trigonometric models, such as the six vertex model, 
or equivalently the XXZ spin chain model (see the references in \cite{JM}), 
we know that a formulation of quantum groups in terms of the Drinfeld 
currents\cite{Drinfeld} provides a convenient framework.
This is because one can construct a free field realization of the 
type I and II intertwining operators starting from a 
free  field realization of the Drinfeld currents.
In addition, the Drinfeld currents have a formal, but deep, resemblance 
to the currents in affine Kac-Moody algebras so that we can easily  
compare the results with those in  conformal field theory.
Hence to perform an algebraic analysis of the elliptic lattice models, it is 
an important step to find a new realization of the both elliptic 
quantum groups $\Aqp(\slnh)$ and $\Bqla(\g)$, $\g$ being an affine 
Lie algebra, in terms of the Drinfeld currents.

In \cite{Konno}, one of the authors has introduced an elliptic analogue of 
the Drinfeld currents of $U_q(\slth)$ independently from the
formulation of the elliptic
 quantum groups. The algebra of the currents is 
called the elliptic algebra
$U_{q,p}(\slth)$. 
Later in \cite{JKOS2}, it has been 
shown that $U_{q,p}(\slth)$ can be regarded  essentially as the Drinfeld 
currents which gives a new realization of the face type elliptic algebra 
 $\Bqla(\slth)$. 
According to this result, the type I and type II vertex operators of 
$U_{q,p}(\slth)$, 
the analogues of the intertwining operators of $\Bqla(\slth)$, 
 have been
realized by the free bosonic fields. The resultant expressions 
coincide with those of the vertex operators of the ABF model 
obtained by Lukyanov and Pugai\cite{LP}. 
Hence a representation theoretical 
foundation to  Lukyanov and Pugai's free field approach to 
the ABF model has been established.   
 
The purpose of this paper is to extend this result to the higher rank case.
We investigate a higher rank elliptic algebra $U_{q,p}(\slnh)$,
and show that $U_{q,p}(\slnh)$ provides a new realization 
of the the face type elliptic algebra  $\Bqla(\slnh)$ in terms of the 
elliptic Drinfeld currents. 

Our strategy is parallel to the one in \cite{JKOS2}. We first give a 
definition of $U_{q,p}(\slnh)$ introducing the new 
currents $K_j(v)\ (1\leq j\leq N)$ (Section \ref{secuqpslnh}).
This gives a completion of the definition of $U_{q,p}(\slnh)$ 
given in Appendix A of \cite{JKOS2}. As an example, a realization 
of $U_{q,p}(\slnh)$ as a certain tensor product of the algebra $U_{q}(\slnh)$ 
and a Heisenberg algebra $\C\{\hat{\H}\}$ is given.
Then we define the ``half currents'' of 
the generating functions (total currents) $E_j(v),\ F_j(v),\ K_j(v)$  
of the algebra $U_{q,p}(\slnh)$ (Section \ref{halfcurrents}). 
The half currents allows us to  
construct a $L$-operator as a Gauss decomposed form of an operator valued
matrix \eqref{def:lhat}. We then argue that 
the thus obtained 
$L$-operator satisfies the $RLL$-relation which characterizes 
the algebra $\Bqla(\slnh)$, when the generators of the mentioned Heisenberg
 algebra are 
reduced to a set of parameters (dynamical parameters) by properly removing  
half of the conjugate variables (Section 5). 
Hence, one can regard the algebra $U_{q,p}(\slnh)$ as a  
tensor product of the algebra $\Bqla(\slnh)$ and the Heisenberg 
algebra $\C\{\hat{\H}\}$.

The $L$-operator and the coalgebra structure of $\Bqla(\slnh)$
allows us to construct a free field realization of the vertex operators of 
$U_{q,p}(\slnh)$, which are extension of the type I and II 
intertwining operators of $\Bqla(\slnh)$  by adding elements of the 
Heisenberg algebra, acting on 
the $U_{q,p}(\slnh)$-modules. 
In the level-one representation, we derived such a realization 
starting from a free field realization of the total currents of 
$U_{q,p}(\slnh)$.  The resultant expressions 
coincide with those of the type I and II vertex operators 
obtained in \cite{AJMP} and \cite{FKQ}. We also show that they 
satisfy the required commutation relations. 
We thus give a representation theoretical meaning to 
the vertex operators of the $A_{N-1}^{(1)}$ type face model\cite{JMO87}. 
Conversely, as a composition of the type I and type II intertwiners, one
can construct a $L$-operator which satisfies the $RLL$-relations
\cite{Miki,FIJKMY}.
As a check of our free field realization, we investigate a connection between 
the two $L$-operators, the one constructed by a composition of the vertex 
operators and the other by the half currents, 
in the level-one representation. We then give a proof of our argument 
in Section 5 at $c=1$. 
 
The article is organized as follows.
In the next section, we review some basic facts on the face type elliptic 
quantum group $\Bqla(\slnh)$. 
In Section \ref{secuqpslnh}, we present a definition of the elliptic algebra 
$U_{q,p}(\slnh)$. New currents $K_j(u)\ (j=1,2,..,N)$ are introduced there. 
A realization of $U_{q,p}(\slnh)$ using the Drinfeld
currents of $U_{q}(\slnh)$ and a Heisenberg algebra
is also given.
In Section \ref{halfcurrents}, we introduce a set of half currents defined 
from $U_{q,p}(\slnh)$ and derive their commutation relations. 
In Section \ref{loperator}, constructing a $L$-operator 
in terms of the half currents, we show that it satisfies the
required $RLL$-relation for $\Bqla(\slnh)$. According to this result, 
in Section \ref{vertexoperators}, 
we discuss a free field realization of the two types of 
vertex operators of the level one $U_{q,p}(\slnh)$-modules. 
In addition, we have four appendices. Appendix A is devoted to a list of
operator product expansions used in the text. In Appendix B, we give a 
proof of some formulae of commutation relations of the half currents.
In Appendix C, we give a derivation of some formulae contained in the
$RLL$-relation. Finally, in Appendix D, we give a summary of the $N$ 
dimensional evaluation representation of $U_{q,p}(\slnh)$.

\section{The Elliptic Quantum Group $\Bqla(\slnhbig)$}\lb{bqlaslnh}
In this section, we give a review on the 
face type elliptic quantum group $\Bqla(\slnh)$ based 
on the results in \cite{JKOS1}. 

\subsection{Notations}
Through this article, we fix a complex number $q\neq0, |q|<1$.
We often use the parameters
\begin{eqnarray*}
p=q^{2r}=e^{-\frac{2\pi i}{\tau}},~~p^*=pq^{-2c}=q^{2r^*}=e^{-\frac{2\pi i}{\tau^*}}~~(r^*=r-c;~ r,r^* \in 
{\mathbb{R}}_{>0},\ r\tau=r^*\tau^*).
\end{eqnarray*}
The following notation is standard:
\begin{eqnarray*}
&&\Theta_p(z)=(z,p)_\infty (pz^{-1};p)_\infty
(p;p)_\infty,\\
&&(z;t_1,\cdots,t_k)_\infty=
\prod_{n_1,\cdots,n_k \geq 0}(1-zt_1^{n_1}\cdots t_k^{n_k}).
\end{eqnarray*}
We also use the Jacobi theta functions
\begin{eqnarray*}
[v]=q^{\frac{v^2}{r}-v}
\frac{\Theta_p(q^{2v})}{(p;p)_\infty^3},~~
[v]^*=q^{\frac{v^2}{r^*}-v}
\frac{\Theta_{p^*}(q^{2v})}{(p^*;p^*)_\infty^3},
\end{eqnarray*}
which satisfy $[-v]=-[v]$ and the quasi-periodicity property
\begin{eqnarray}
&&~[v+r]=-[v],~~[v+r\tau]=-e^{-\pi i \tau-\frac{2\pi i v}{r}}[v].
\end{eqnarray}
We take the normalization of the theta function to be 
\begin{eqnarray}
\oint_{C_0}\frac{dz}{2\pi i z}\frac{1}{[-v]}=1,
\end{eqnarray}
where $C_0$ is a simple closed curve in the $v$-plane encircling
$v=0$ anticlockwise.
The same holds for $[v]^*$, with $r$
replaced by $r^*$, except for the normalization
\be
&&\oint_{C_0}\frac{dz}{2\pi i z}\frac{1}{[-v]^*}=\frac{[v]}{[v]^*}
\Bigl\vert_{v\to 0}.
\en

\subsection{Definition of the elliptic quantum group $\Bqla(\slnh)$}
Let $U_q=U_{q}(\slnh)$ be the standard affine quantum group. Namely, $U_q(\slnh)$ is a quasi-triangular
Hopf algebra equipped with the standard coproduct $\Delta$,  counit $\vep$, 
antipode $S$ and universal $R$ matrix $\cR$. 
Our conventions on the coalgebra structure follows \cite{JKOS1}. 
Let $\h$ and $\bar{\h}$ be the Cartan subalgebras of $\slnh$ and $\sln$, 
respectively.
We denote a basis and its dual basis of ${\h}$ by $\{\hat{h}_l\}$ and $\{\hat{h}^l\}$, respectively. More explicitly, they are given by $\{\hat{h}_l\}=\{d, c, h_j \}$ and $\{\hat{h}^l\}=\{c, d, h^j \}\ (1\leq j\leq N-1)$, where $c$ and 
$d$ are a central element and a derivation operator of $\slnh$, 
respectively, and 
$\{{h}_j\}$ and $\{{h}^j\}$ are a basis  and a dual basis of 
$\bar{\h}$. 

The face type elliptic quantum group $\Bqla(\slnh)$ is a quasi-Hopf
deformation of $U_{q}(\slnh)$ by the face type twistor $F(\la)\ (\la \in \h)$.
The twistor $F(\la)$ is an invertible element in $U_q\otimes U_q$ 
satisfying 
\bea
&&({\rm id}\otimes \vep )F(\la)=1=F(\la)(\vep \otimes {\rm id}),\\
&&F^{(12)}(\la)(\Delta\otimes {\rm id})F(\la)
=F^{(23)}(\la+h^{(1)})({\rm id} \otimes \Delta)F(\la).
\lb{facecocy}
\lb{epF}
\ena
where $\la=\sum_l\la_l\hat{h}^l\  (\la_l\in \C)$, 
$\la+h^{(1)}=\sum_l(\la_l+\hat{h}_l^{(1)})\hat{h}^l$ and 
$\hat{h}_l^{(1)}=\hat{h}_l\otimes 1\otimes 1$.  
An explicit construction of the twistor $F(\la)$ is given in \cite{JKOS1}.
A quasi-Hopf deformation means that as an associative algebra,  
$\Bqla(\slnh)$ is isomorphic to $U_q(\slnh)$, 
but the coalgebra structure is deformed. Namely, the coproduct is changed to
the new one given by
\bea
&&\Delta_\la(x)=F(\la)\Delta(x)F(\la)^{-1}\qquad \forall x \in U_q(\slnh).
\ena
$\Delta_\la$ satisfies a weaker coassociativity
\bea
&&({\rm id}\otimes \Delta_\la )\Delta_\la(x)=\Phi(\la)
(\Delta_\la \otimes {\rm id})\Delta_\la(x)\Phi(\la)^{-1}\qquad \forall\ x \in U_q(\slnh),\\
&&\Phi(\la)=F^{(23)}(\la)F^{(23)}(\la+h^{(1)})^{-1}.
\ena
The universal $R$-matrix is also deformed to
\bea
&&\cR(\la)=F^{(21)}(\la)\cR F^{(12)}(\la)^{-1}.
\ena

\begin{dfn}{\bf (Elliptic quantum group $\Bqla(\slnh)$)}\cite{JKOS1}~~
The face type elliptic quantum group $\Bqla(\slnh)$ is a quasi-triangular
quasi-Hopf algebra 
$(\Bqla(\slnh),$
$~\Delta_\la, \vep,~S,~\Phi(\la),~\alpha,~\beta,~\cR(\la))$,
where $\alpha,\ \beta$ are defined by 
\bea
&&\alpha=\sum_iS(k_i)l_i,\quad \beta=\sum_i m_iS(n_i).
\ena
Here we set $\sum_ik_i\otimes l_i=F(\la)^{-1},\ 
\sum_im_i\otimes n_i=F(\la)$.
\end{dfn}

A characteristic feature of $\Bqla(\slnh)$ is that the 
universal $R$ 
matrix $\cR(\la)$ satisfies the dynamical Yang-Baxter equation.
\begin{equation}
\cR^{(12)}(\la+ h^{(3)} )\cR^{(13)}(\la)\cR^{(23)}(\la+h^{(1)})
=\cR^{(23)}(\la)\cR^{(13)}(\la+h^{(2)})\cR^{(12)}(\la).
\label{DYBE}
\end{equation}
Let $(\pi_{V,z},V_z),\ V_z=V\otimes \C[z,z^{-1}]$ be a (finite dimensional) 
evaluation representation 
of $U_q$. Taking images of $\cR$, we have a $R$-matrix 
$R^+_{VW}(z,\la)$ and a $L$-operator $L^+_V(z,\la)$ as follows.
\bea
&&R^{+}_{VW}(z_1/z_2,\la)=\left(\pi_{V,z_1}\otimes\pi_{W,z_2}\right)
q^{c\otimes d+d\otimes c}\cR(\la),\\
&&L_V^+(z,\la)=\left(\pi_{V,z}\otimes {\rm id} \right)q^{c\otimes d+d\otimes c}\cR(\la).
\ena
Then from \eqref{DYBE}, we have the following dynamical $RLL$-relation. 
\bea
&&R^{+}_{VW}(z_1/z_2,\la+ h)
L_V^{+}(z_1,\la)L_W^{+}(z_2,\la+ h^{(1)})
=
L_W^{+}(z_2,\la)L_V^{+}(z_1,\la+ h^{(2)})
R^{+}_{VW}(z_1/z_2,\la).\nn\\
\lb{DRLL}
\ena
Note that in $\Bqla(\slnh)$, 
$L^+_V(z,\la)$ and $L^-_V(z,\la)=\left(\pi_{V,z}\otimes {\rm id} \right)
\cR^{(21)}(\la)^{-1}q^{-c\otimes d-d\otimes c}$ are not independent 
operators (Proposition 4.3 in \cite{JKOS1}). Hence just one 
dynamical $RLL$-relation \eqref{DRLL} characterizes the algebra 
$\Bqla(\slnh)$ completely in the sense of Reshetikhin and Semenov-Tian-Shansky
\cite{RS}.

Hereafter we parametrize the dynamical variable $\lambda$ as
\bea 
&&\lambda=(r^*+N)d+s'c+\sum_{j=1}^{N-1}(s_j+1)h^{j}
\qquad 
(s'\in \C,\ r^*\equiv r-c).\lb{sj}
\ena 
Under this, we set  $F(r^*,\{s_j\})\equiv F(\la)$ and 
$\cR(r^*,\{s_j\}) \equiv \cR(\la)$.
Since $c$ is central, no $s'$ dependence should appear. 
The dynamical shift $\lambda \to \lambda+h$  
with $h=cd+\sum_{j=1}^Nh_jh^j$, 
changes  the universal $R$-matrix $\cR(r^*,\{s_j\})$ to 
$\cR(r,\{s_j+h_{j}\})\equiv \cR(\lambda+h)$. Note $r^*=r-c$.

Let us now take $(\pi_{V,z}, V_z)$ to 
be the  evaluation representation associated with the 
vector representation $V\cong \C^N$ of $U_q(\sln)$ 
(see Appendix \ref{Evaluation}).
We set 
\be
&&R^+(v,s+h)=(\pi_{V,z_1}\otimes \pi_{V,z_2})q^{c\otimes d+d\otimes c}\cR(r,\{s_j+h_j\}),\\
&&L^+(v,s)=(\pi_{V,z}\otimes \id)q^{c\otimes d+d\otimes c}\cR(r^*,\{s_j\}),
\en
where $z_i=q^{2v_i}\ (i=1,2)$, $v=v_1-v_2$. 
One can obtain the finite dimensional representation of the twistor
$F(r,\{s_j\})$ 
by solving the difference equation for $(\pi_{V,z_1}\otimes \pi_{V,z_2})F(r,\{s_j+h_j\})$ (Eq.(2.30) in \cite{JKOS1}) 
derived by using the explicit 
realization of $F(\la)$, under the parametrization \eqref{sj}.
Then noting the relation $\cR(r,\{s_j+h_j\})=F^{(21)}(r,\{s_j+h_j\})\cR F^{(12)}(r,\{s_j+h_j\})^{-1}$,
we obtain the $R$-matrix $R^{+}(v,s+h)$, up to a certain gauge transformation,
as  
\begin{eqnarray}
R^{+}(v,s+h)&=&\rho^{+}(v)\bar{R}(v,s+h),\lb{rmatfull}\\
\bar{R}(v,s+h)&=&
\sum_{j=1}^{N}
E_{jj}\otimes E_{jj}+
\sum_{1 \leq j<l \leq N}
\left(b_{}(v,s_{j,l}+h_{j,l})
E_{jj}
\otimes E_{ll}+
\bar{b}_{}(v)
E_{ll}\otimes E_{jj}
\right)\nonumber\\
&&+
\sum_{1 \leq j<l \leq N}
\left(c_{}(v,s_{j,l}+h_{j,l})
E_{jl}\otimes E_{lj}+
\bar{c}_{}
(v,s_{j,l}+h_{j,l})E_{lj}\otimes E_{jl}
\right),\lb{rmat}
\end{eqnarray}
where $s_{j,l}=\sum_{m=j}^{l-1}s_j,\ h_{j,l}=\sum_{m=j}^{l-1}h_j \ (1\leq j<l \leq N)$ and 
\begin{eqnarray}
&&b_{}(u,s)=
\frac{[s+1]
[s-1][u]
}{
[s]^{2} [u+1]
},
~~~\bar{b}_{}(u)=\frac{[u]}{[u+1]},\\
&&c_{}(u,s)=
\frac{[1][s+u]}{
[s][u+1]},
~~~\bar{c}_{}(u,s)=
\frac{[1][s-u]}{
[s][u+1]}.\lb{rmatcomp}
\end{eqnarray}
The function $\rho^+(v)$ is chosen as 
\bea
&&\rho^+(v)=q^{\frac{N-1}{N}}
z^{\frac{N-1}{rN}}
\frac{\{pq^2z\}
\{pq^{2N-2}z\}
\{1/z\} \{q^{2N}/z\}
}{
\{pz\}
\{pq^{2N}z\}
\{q^2/z\} \{q^{2N-2}/z\}
},\label{def:rho}
\ena
where 
\begin{eqnarray}
\{z\}=(z;p,q^{2N})_\infty.
\end{eqnarray}
Up to a gauge transformation, the $R$-martrix $R^+(v,s+h)$ is nothing but the Boltzmann 
weight of the $A_{N-1}^{(1)}$ type face model  
introduced in \cite{JMO87}.
The $R$-matrix  $R^{+*}(v,s)=(\pi_{V,z_1}\otimes \pi_{V,z_2})\cR(r^*,\{s_j\})$
 is obtained from $R^{+}(v,s)$ by the 
replacements $r \to r^*$.
Hence, under the parametrization \eqref{sj},  
the dynamical $RLL$-relation takes the form 
\bea
&&R^{+(12)}(v,s+h)
L^{+(1)}(v_1,s)L^{+(2)}(v_2,s+ h^{(1)})
=
L^{+(2)}(v_2,s)L^{+(1)}(v_1,s+h^{(2)})
R^{+*(12)}(v,s).\nn\\
\lb{DRLL2}
\ena

\subsection{Intertwining operators}

Let $\F, \F'$ be highest weight $U_q$-modules.
We denote the type-I and type II intertwining operators of $U_q$-modules 
by $\Phi(z)$
and 
$\Psi^*(z)$, respectively. 
\begin{eqnarray}
\Phi(z):~{\cal F} \longrightarrow {\cal F}'\otimes W_z,~~\qquad
\Psi^*(z):~W_z \otimes {\cal F} \longrightarrow
{\cal F}'.
\end{eqnarray}
Twisting these operators by $F(r^*,s)$,
we obtain the corresponding intertwining operators ${\Phi}(v,s)$ 
and ${\Psi}^*(v,s)$ of $\Bqla$-modules.
\begin{eqnarray}
&&{\Phi}_W(v,s)=(\id \otimes \pi_{W,z})F(r^*,\{s_j\}) \Phi(z),\lb{intI}\\
&&{\Psi}_W^*(v,s)=\Psi^*(z) (\pi_{W,z} \otimes \id)F(r^*,\{s_j\})^{-1}.
\lb{intII}
\end{eqnarray}
From the intertwining relation satisfied by $\Phi(z)$ and $\Psi^*(z)$, 
one can derive  
the following dynamical intertwining relation for the  new intertwiners 
\cite{JKOS1}.
\begin{eqnarray}
&&{\Phi}_W^{(3)}(v_2+\frac{c}{2},s)
L_V^{+(1)}(v_1,s)=R^{+(13)}_{VW}(v,s+h)
L_V^{+(1)}(v_1,s){\Phi}_W^{(3)}(v_2+\frac{c}{2},s+h^{(1)}),
\nonumber\\
\lb{dintrelI}\\
&&L_V^{+(1)}(v_1,s){\Psi}_W^{*(2)}(z_2,s+h^{(1)})=
{\Psi}_W^{*(2)}(z_2,s)L_V^{+(1)}
(v_1,s+h^{(2)})R_{VW}^{+*(12)}(v_1-v_2,s).
\nonumber
\\
\lb{dintrelII}
\end{eqnarray}
Note that \eqref{dintrelI} and \eqref{dintrelII} are the relations for
the operators $V_{z_1}\otimes {\cal F} \to 
V_{z_1}\otimes {\cal F} \otimes W_{z_2}$ and 
$V_{z_1}\otimes W_{z_2}\otimes {\cal F} \to 
V_{z_1}\otimes {\cal F}$, respectively. 

\section{The Elliptic Algebra $U_{q,p}(\slnhbig)$}\lb{secuqpslnh}
In this section,  we give a definition of the elliptic algebra $U_{q,p}(\slnh)$. To define the algebra, we  follows mainly the idea given in  Appendix A of 
\cite{JKOS2}.
Namely, we first introduce the elliptic currents $e_i(z,p),\ f_i(z,p)$ and 
$\psi^{\pm}_i(z,p)$ of $U_q(\slnh)$ by modifying the Drinfeld currents of 
$U_q(\slnh)$. 
Then we extend them to the currents of $U_{q,p}(\slnh)$ by taking a 
tensor product with a Heisenberg algebra $\C\{\hat{\H}\}$ given
 in Section 3.4.1. 
Our definition is an extended version of the one given in \cite{JKOS2}
 introducing new currents $K_j(v)\ (1\leq j\leq N)$. The currents 
 $\{ K_j(v) \}$ play an essential role in 
the construction of the $L$-operators (Section 5).

\subsection{Drinfeld currents of $U_q(\slnh)$}\lb{drinfelduq}
Let us first recall the Drinfeld currents of
$U_q(\widehat{\goth{sl}}_N)$ \cite{Drinfeld}.
We use the standard symbol of $q$-integer
\begin{eqnarray}
[n]_q=\frac{q^n-q^{-n}}{q-q^{-1}}.
\end{eqnarray}
We also use the symbol $A=(A_{j k})$ to express the Cartan matrix of $\sln$.
\begin{df}{\bf (Drinfeld currents)}
~~~~The algebra $U_q(\widehat{\goth{sl}}_N)$
is a $\mathbb{C}$-algebra generated by the generators
$h_i,\  a_{i,m},\ x_{i,n}^\pm\ (i=1,\cdots,N-1: m 
\in {{\Z}_{\neq 0}},\ n \in {\mathbb{Z}}),\ c,\ d$.
In terms of the generating functions

\begin{eqnarray}
&&x_i^\pm(z)=\sum_{n \in {\mathbb{Z}}}x_{i,n}^\pm z^{-n},\\
&&\psi_i(q^{\frac{c}{2}}z)=q^{h_i}
\exp\left(
(q-q^{-1})\sum_{m>0}a_{i,m}z^{-m}\right),\\
&&\varphi_i(q^{-\frac{c}{2}}z)=q^{-h_i}
\exp\left(-
(q-q^{-1})\sum_{m>0}a_{i,-m}z^{m}\right)\qquad (i=1,\cdots,N-1),
\end{eqnarray}
the defining relations of $U_q(\slnh)$ are given by
\begin{eqnarray}
&&~{c}: {\rm central},\label{def:D1}\\
&&~[h_i,d]=[d, h_i]=[d, a_{i,m}]=[a_{i,m}, d]=0,\label{def:D2}\\
&&~[d, x_{i,n}^\pm]=n\ x_{i,n}^\pm,
~[h_i, a_{j,m}]=[a_{j,m}, h_i]=0,\label{def:D3}\\
&&~[h_i, x_j^\pm(z)]=\pm A_{ij}\ x_j^\pm(z),\label{def:D4}\\
&&~[a_{i,m},a_{j,n}]=\frac{[A_{ij}m]_q
[c m]_q}{m}\ q^{-c|m|}
\delta_{n+m,0},\label{def:D5}\\
&&~[a_{i,m},x_j^+(z)]=\frac{[A_{ij}m]_q}{m}\ q^{-c|m|}z^m
x_j^+(z),\label{def:D6}\\
&&~[a_{i,m},x_j^-(z)]=-\frac{[A_{ij}m]_q}{m}\ z^m x_j^-(z),
\label{def:D7}\\
&&~(z_1-q^{\pm A_{ij}}z_2)x_i^\pm(z_1)x_j^\pm(z_2)=
(q^{\pm A_{ij}}z_1-z_2)x_j^\pm(z_2)x_i^\pm(z_1),
\label{def:D8}\\
&&~[x_i^+(z_1),x_j^-(z_2)]=\frac{\delta_{i,j}}{q-q^{-1}}
\left(\delta(q^{-c}z_1/z_2)\psi_i(q^{\frac{c}{2}}z_2)
-\delta(q^c z_1/z_2)\varphi_i(q^{-\frac{c}{2}}z_2)\right),
\label{def:D9}\\
&&~(x_i^\pm(z_1)x_i^\pm(z_2)x_j^\pm(z)-
[2]_q\ x_i^\pm(z_1)x_j^\pm(z)x_i^\pm(z_2)+
x_j^\pm(z)x_i^\pm(z_1)x_i^\pm(z_2))\nonumber\\
&&~+
(x_i^\pm(z_2)x_i^\pm(z_1)x_j^\pm(z)-
[2]_q\ x_i^\pm(z_2)x_j^\pm(z)x_i^\pm(z_1)+
x_j^\pm(z)x_i^\pm(z_2)x_i^\pm(z_1))=0,\nonumber
\\
&&~~~~~~~~~~~~~~~~~~~~~~~~~~~~
{\rm for}~~~|i-j|=1.\label{def:D10}
\end{eqnarray}
Here $\delta(z)$ denotes the delta function
$\delta(z)=\sum_{m \in {\mathbb{Z}}}z^m$.
We call the generators ${h_j},\ a_{j,m},\ x_{j,n}^\pm,\  c,\ d$ 
the Drinfeld generators of $U_q(\widehat{\goth{sl}}_N)$ 
and the generating functions $x_i^\pm(z), \psi_i(z)$ and $\varphi_i(z)$
the Drinfeld currents.
\end{df}

\subsection{Elliptic currents of $U_q(\widehat{\goth{sl}}_N)$}
We next introduce an elliptic modification of the currents
$x_i^\pm(z), \psi_i(z)$ and $\varphi_i(z)$ according to \cite{JKOS2}.

Let us define the auxiliary currents
$u_i^\pm(z,p)$ by
\begin{eqnarray}
&&u_i^+(z,p)=
\exp\left(
\sum_{m>0}\frac{1}{[r^* m]_q}a_{i,-m}(q^r z)^m\right),
\label{def:u1}\\
&&u_i^-(z,p)=
\exp\left(
-\sum_{m>0}\frac{1}{[r m]_q}a_{i,m}(q^{-r}z)^{-m}\right).
\label{def:u2}
\end{eqnarray}
\begin{prop}~~~
The following commutation relations hold.
\begin{eqnarray}
&&u_i^+(z_1,p)x_j^+(z_2)=
\frac{(p^*q^{A_{ij}}z_1/z_2;p^*)_\infty}
{(p^*q^{-A_{ij}}z_1/z_2;p^*)_\infty}
x_j^+(z_2)u_i^+(z_1,p),\\
&&u_i^+(z_1,p)x_j^-(z_2)=
\frac{(p^*q^{c-A_{ij}}z_1/z_2;p^*)_\infty}
{(p^*q^{c+A_{ij}}z_1/z_2;p^*)_\infty}x_j^-(z_2)u_i^+(z_1,p),
\\
&&u_i^-(z_1,p)x_j^+(z_2)=
\frac{(pq^{-c-A_{ij}}z_2/z_1;p)_\infty}
{(pq^{-c+A_{ij}}z_2/z_1;p)_\infty}
x_j^+(z_2)u_i^+(z_1,p),\\
&&u_i^-(z_1,p)x_j^-(z_2)=
\frac{(pq^{A_{ij}}z_2/z_1;p)_\infty}
{(pq^{-A_{ij}}z_2/z_1;p)_\infty}x_j^-(z_2)u_i^-(z_1,p),
\\
&&u_i^+(z_1,p)u_j^-(z_2,p)=\frac{(pq^{-c-A_{ij}}z_1/z_2;p)_\infty}{(pq^{-c+A_{ij}}z_1/z_2;p)_\infty}
\frac{(p^*q^{c+A_{ij}}z_1/z_2;p^*)_\infty}
{(p^*q^{c-A_{ij}}z_1/z_2;p^*)_\infty}
u_j^-(z_2,p)u_i^+(z_1,p).\nonumber\\
\end{eqnarray}
\end{prop}

\begin{df}{\bf (Elliptic currents)}~~~~
Let us define ``dressed'' currents
${e}_i(z,p), {f}_i(z,p)$, 
${\psi}_i^\pm(z,p),$
$(i=1,\cdots,N-1)$ by
\begin{eqnarray}
&&e_i(z,p)=u_i^+(z,p)x_i^+(z),\label{def:e}\\
&&f_i(z,p)=x_i^-(z)u_i^-(z,p),\label{def:f}\\
&&\psi_i^+(z,p)=u_i^+(q^{\frac{c}{2}}z,p)\psi_i(z)
u_i^-(q^{-\frac{c}{2}}z,p)
,\label{def:psi+}\\
&&\psi_i^-(z,p)=
u_i^+(q^{-\frac{c}{2}}z,p)\varphi_i(z)
u_i^-(q^{\frac{c}{2}}z,p).\label{def:psi-}
\end{eqnarray}
We call the currents 
${e}_i(z,p), {f}_i(z,p)$ and $ {\psi}_i^\pm(z,p)$
the elliptic currents of $U_q(\widehat{\goth{sl}}_N)$.
\end{df}
The reason why we call ``elliptic'' is because the dressing 
operation specified by $u^{\pm}_i(z,p)$ changes the commutation relation of the 
Drinfeld currents to the elliptic ones.
\begin{prop}\lb{ellCR}~~~The elliptic currents satisfy the following
relations.
\begin{eqnarray}
&&z_1 \Theta_{p^*}(q^{A_{ij}}z_2/z_1)e_i(z_1,p)e_j(z_2,p)=
-z_2 \Theta_{p^*}(q^{A_{ij}}z_1/z_2)e_j(z_2,p)e_i(z_1,p),\lb{ee}\\
&&\nn\\
&&
z_1 \Theta_{p}(q^{A_{ij}}z_2/z_1)f_i(z_1,p)f_j(z_2,p)=
-z_2 \Theta_{p}(q^{A_{ij}}z_1/z_2)f_j(z_2,p)f_i(z_1,p),\lb{ff}\\
&&\nn\\
&&[e_i(z_1,p),f_j(z_2,p)]=\frac{\delta_{i,j}}{q-q^{-1}}
\left(\delta(q^{-c}z_1/z_2)
\psi_j^+(q^{\frac{c}{2}}z_2,p)-
\delta(q^cz_1/z_2)
\psi_j^-(q^{-\frac{c}{2}}z_2,p)
\right),\lb{ef}\\&&
\nonumber\\
&&q^{-h_j}\psi_j^+(q^{-r+\frac{c}{2}}z,p)=
q^{h_j}\psi_j^-(q^{r-\frac{c}{2}}z,p),\lb{psippsim}\\
&&\nn\\
&&\psi_i^+(z_1,p)\psi_j^+(z_2,p)=
\frac{\Theta_p(q^{-A_{ij}}z_1/z_2)\Theta_{p^*}(q^{A_{ij}}z_1/z_2)}
{\Theta_{p}(q^{A_{ij}}z_1/z_2)\Theta_{p^*}(q^{-A_{ij}}z_1/z_2)}
\psi_j^+(z_2,p)\psi_i^+(z_1,p)
,\lb{psipsi}\\
&&\nn\\
&&\psi_i^+(z_1,p)e_j(z_2,p)=
\frac{\Theta_{p^*}(q^{A_{ij}+\frac{c}{2}}z_1/z_2)}{
\Theta_{p^*}(q^{-A_{ij}+\frac{c}{2}}z_1/z_2)}
e_j(z_2,p)\psi_i^+(z_1,p)
,\lb{psie}\\
&&\nn\\
&&\psi_i^+(z_1,p)f_j(z_2,p)=
\frac{\Theta_{p}(q^{-A_{ij}-\frac{c}{2}}z_1/z_2)}{
\Theta_{p}(q^{A_{ij}-\frac{c}{2}}z_1/z_2)}
f_j(z_2,p)\psi_i^+(z_1,p)
,\lb{psif}\\
&&\nn\\
&&\frac{(p^*q^2z_2/z_1;p^*)_\infty}
{(p^*q^{-2}z_2/z_1;p^*)_\infty
} \left\{
\frac{(p^*q^{-1}z/z_1;p^*)_\infty 
(p^*q^{-1}z/z_2;p^*)_\infty}
{(p^*qz/z_1;p^*)_\infty 
(p^*qz/z_2;p^*)_\infty}e_i(z_1,p)
e_i(z_2,p)
e_j(z,p)\right.\nonumber\\
&&\qquad\qquad -\left.[2]_q\frac{(p^*q^{-1}z/z_1;p^*)_\infty 
(p^*q^{-1}z_2/z;p^*)_\infty}
{(p^*qz/z_1;p^*)_\infty 
(p^*qz_2/z;p^*)_\infty}
e_i(z_{1},p)e_j(z,p)e_i(z_{2},p)
\right.\nonumber\\
&&\qquad\qquad +\left.
\frac{(p^*q^{-1}z_1/z;p^*)_\infty 
(p^*q^{-1}z_2/z;p^*)_\infty}
{(p^*qz_1/z;p^*)_\infty 
(p^*qz_2/z;p^*)_\infty}
e_j(z,p)e_i(z_{1},p)
e_i(z_{2},p)
\right\}+(z_1 \leftrightarrow z_2)=0,\nn\\
&&\lb{serre}\\
&&\frac{(pq^{-2}z_2/z_1;p)_\infty}
{(pq^{2}z_2/z_1;p)_\infty
} \left\{
\frac{(pq z/z_1;p)_\infty 
(pq z/z_2;p)_\infty}
{(pq^{-1} z/z_1;p)_\infty 
(pq^{-1} z/z_2;p)_\infty}
f_i(z_1,p)f_i(z_2,p)f_j(z,p)\right.\nonumber\\
&&\qquad\qquad
-\left.[2]_q\frac{(pq z/z_1;p)_\infty 
(pq z_2/z;p)_\infty}
{(pq^{-1} z/z_1;p)_\infty 
(pq^{-1} z_2/z;p)_\infty}
f_i(z_{1},p)f_j(z,p)f_i(z_{2},p)
\right.\nonumber
\\
&&\qquad\qquad+\left.
\frac{(pq z_1/z;p)_\infty 
(pq z_2/z;p)_\infty}
{(pq^{-1} z_1/z;p)_\infty 
(pq^{-1} z_2/z;p)_\infty}
f_j(z,p)f_i(z_{1},p)
f_i(z_{2},p)
\right\} +(z_1 \leftrightarrow z_2)=0,\nn\\
&&~~~~~~~~~~~\qquad\qquad{\rm for}~~~|i-j|=1.\lb{serrf}
\end{eqnarray}
\end{prop}

\subsection{New currents $k_j(v)$}
In this subsection, 
we consider a decomposition of the elliptic currents
$\psi^\pm_j(z,p)\ (1\leq j\leq N-1)$ corresponding to the decomposition
 \eqref{root}. For this purpose, we introduce new currents 
 $k_j(v)\ (1\leq j\leq N)$.

We first note that the currents $\psi^\pm_j(z,p)$ are 
expressed by using the Drinfeld generators $a_{j,m}$ as follows. 
\bea
&&\psi_j^\pm(q^{\mp (r -\frac{c}{2})}z,p)=q^{\pm h_j}:\exp\left\{
-\sum_{m\neq0}\frac{1}{[r^*m]_q}
b_{j,m}
(q^{N-j}z)^{-m}
\right\}:,\lb{psib0}
\ena
where we set
\begin{eqnarray}
b_{j,m}=\left\{
\begin{array}{cc}
\frac{[r^*m]_q}{[rm]_q}a_{j,m}& m>0\\
q^{c|m|}a_{j,m}& m<0.
\end{array}
\right. \label{def:Boson2}
\end{eqnarray}
The colons in \eqref{psib0} denote the standard normal ordering.

Let us introduce new generators,
$B_m^j\ (j=1,\cdots,N; m \in {\mathbb{Z}})$,
according to the formula
\begin{eqnarray}
-B_m^j+B_m^{j+1}=\frac{m}{[m]_q}b_{j,m} q^{(N-j)m},~~\qquad
\sum_{j=1}^N q^{2jm} B_m^j=0,\label{def:Boson1}
\end{eqnarray}
or more explicitly, 
\begin{eqnarray}
B_m^j=\frac{m}{[m]_q[Nm]_q}\left(
\sum_{k=1}^{j-1}[km]_q\ b_{k,m}
-q^{Nm} \sum_{k=j}^{N-1}[(N-k)m]_q\ b_{k,m}
\right).\label{def:Boson3}
\end{eqnarray}
From this and \eqref{def:D5}-\eqref{def:D7}, 
we derive the following commutation relations.
\begin{prop}~~~
For $m,m' \in {{\Z}_{\neq 0}},\ j,k=1,\cdots,N$,
the following commutation relations hold.
\begin{eqnarray}
&&~[B_m^j,B_{m'}^k]=m \delta_{m+m',0}
\frac{[r^*m]_q[cm]_q}{[rm]_q[m]_q[Nm]_q}\times
\left\{\begin{array}{cc}
[(N-1)m]_q &\ (j=k)\\
-q^{-mN{\rm sgn(j-k)}}[m]_q&\ (j\neq k),
\end{array}\right.\\
&&~[B_m^j,x_j^\pm(z)]=\mp
q^{m(N+1-j-{c})}
z^m x_j^\pm(z)\times \left\{
\begin{array}{cc}
\frac{[r^*m]_q}{[rm]_q}&\ (m>0) \\
q^{c m}&\  (m<0),
\end{array}\right.
\\
&&~[B_m^{j+1},x_j^\pm(z)]=\pm q^{m(N-1-j-{c})}z^m x_j^\pm(z)\times\left\{
\begin{array}{cc}
\frac{[r^*m]_q}{[rm]_q}&\ (m>0) \\
q^{c m}&\  (m<0),
\end{array}\right.\\
&&~[B_m^k,x_j^\pm(z)]=0~~\qquad ( k\neq j,j+1).
\end{eqnarray}
\end{prop}

We now define new currents $k_j(z,p)\ (1\leq j\leq N)$  by
\begin{eqnarray}
k_j(z,p)= 
:\exp\left(\sum_{m\neq0}\frac{[m]_q}
{m [r^*m]_q}B_m^j z^{-m}\right):.
\label{def:k}
\end{eqnarray}
Then, from \eqref{psib0} and \eqref{def:Boson1}, 
we have the following decomposition.
\begin{eqnarray}
&&\psi^\pm_j(q^{\pm (r - \frac{c}{2})}z,p)=
{\kappa}\ q^{\pm h_j} k_j(q^{N-j}z,p) 
k_{j+1}(q^{N-j}z,p)^{-1},\nonumber\\
&&\qquad {\kappa}=\frac{(p;p)_\infty (p^*q^2;p^*)_\infty}
{(p^*;p^*)_\infty (pq^2;p)_\infty}.\lb{psikk}
\end{eqnarray} 
It is also easy to verify the following commutation relations. 
\begin{prop}~~~
\begin{eqnarray}
&&k_{j}(z_1,p)k_j(z_2,p)=
\left(\frac{z_1}{z_2}\right)^{\frac{N-1}{N}(\frac{1}{r}-\frac{1}{r^*})}
{{\rho}(v_1-v_2)}
k_j(z_2,p)k_{j}(z_1,p),\lb{kjkj}
\\
&&k_{j_1}(z_1,p)k_{j_2}(z_2,p)=\left(\frac{z_1}{z_2}\right)^{\frac{N-1}{N}(\frac{1}{r}-\frac{1}{r^*})}
{\rho}(v_1-v_2)
\frac{\Theta_p(z_1/z_2)}{\Theta_{p}(q^{-2}z_1/z_2)}
\frac{\Theta_{p^*}(q^{-2}z_1/z_2)}{\Theta_{p^*}(z_1/z_2)}
k_{j_2}(z_2,p)k_{j_1}(z_1,p),\nonumber\\
&&~~~~~~~~~~~~~~~~~~~~\qquad\qquad(1\leqq j_1<j_2 \leqq N),\lb{kj1kj2}\\
&& k_1(z,p)k_2(q^2z,p)\cdots k_{N-1}(q^{2(N-2)}z,p) k_N(q^{2(N-1)}z,p)
={c}_N^{{(N-1)}},\lb{prodk}\\
&&\nn\\
&&k_{j}(z_1,p)e_{j}(z_2,p)=
\frac{\Theta_{p^*}(q^{j-N+r^*}z_1/z_2)}
{\Theta_{p^*}(q^{j-N+r^*-2}z_1/z_2)}
e_{j}(z_2,p)k_{j}(z_1,p),\lb{kjej}
\\&&
k_{j+1}(z_1,p)e_{j}(z_2,p)=
\frac{\Theta_{p^*}(q^{j-N+r^*}z_1/z_2)}
{\Theta_{p^*}(q^{j-N+r^*+2}z_1/z_2)}e_{j}(z_2,p)k_{j+1}(z_1,p),\lb{kjp1ej}
\\&&
k_{i}(z_1,p)e_{j}(z_2,p)=e_{j}(z_2,p)k_{i}(z_1,p),
~\qquad (i\neq j,j+1),\lb{kiej}\\
&&\nn\\
&&
k_{j}(z_1,p)f_{j}(z_2,p)=
\frac{\Theta_p(q^{j-N+r-2}z_1/z_2;p)}
{\Theta_p(q^{j-N+r}z_1/z_2;p)}f_{j}(z_2,p)k_{j}(z_1,p),
\lb{kjfj}\\
&&
k_{j+1}(z_1,p)f_{j}(z_2,p)=
\frac{\Theta_p(q^{j-N+r+2}z_1/z_2)}
{\Theta_p(q^{j-N+r}z_1/z_2)}f_{j}(z_2,p)k_{j+1}(z_1,p),\lb{kjp1fj}
\\&&
k_{i}(z_1,p)f_{j}(z_2,p)=f_{j}(z_2,p)k_{i}(z_1,p),
~\qquad (i\neq j,j+1).\lb{kifj}
\end{eqnarray}
Here we set
\begin{eqnarray}
&&\rho(v)=\frac{\rho^+(v)}{\rho^{+*}(v)},\label{rhofunc}\\
&&{c}_N=\frac{\{pq^{2N+4}\} \{pq^{2N}\} 
\{p^*q^{2N+2}\}^{* 2}}
{\{pq^{2N+2}\}^2 \{p^*q^{2N+4}\}^* 
\{p^*q^{2N}\}^*}
\end{eqnarray}
with $\rho^+(v)$ given in \eqref{def:rho} and  $\rho^{+*}(v)=\rho^+(v)|_{r\to r^*}$.
\end{prop}

\subsection{Definition of the 
elliptic algebra $U_{q,p}(\widehat{\goth{sl}}_N)$}
Now we give a definition of the elliptic algebra
$U_{q,p}(\widehat{\goth{sl}}_N)$ by considering a tensor product of the 
elliptic currents of $U_q(\slnh)$ with 
a Heisenberg algebra.
In order to keep the defining relations of 
the algebra $U_{q,p}(\slnh)$ with the new currents $K_j(v)$ same as those given in Appendix A of 
\cite{JKOS2}, 
we need to make a central extension of the Heisenberg algebra. 

\subsubsection{The Heisenberg algebra $\C\{\H\}$ and its extension 
$\C\{\hat{\H}\}$}\lb{heisenberg}

Let $ \epsilon_j\ (1\leq j\leq N)$
be the orthonormal basis in $\R^N$ with the inner product
$\langle \epsilon_j, \epsilon_k \rangle=\delta_{j,k}$.
Setting 
\bea
&&\bar{\epsilon}_j=\epsilon_j-\epsilon,~\qquad
\epsilon=\frac{1}{N}\sum_{j=1}^N \epsilon_j,
\ena
we have the weight lattice ${P}$ of $A_{N-1}^{(1)}$ 
\begin{eqnarray}
{P}=\oplus_{j=1}^N {\mathbb{Z}}\ \bar{\epsilon}_j.
\end{eqnarray}
Then the simple roots  $\alpha_j\ (1\leq j\leq N-1)$ of $\sln$ are given
by
\begin{eqnarray}
\alpha_j=-\bar{\epsilon}_j+
\bar{\epsilon}_{j+1}.\lb{root}
\end{eqnarray}
Let us introduce operators $h_\alpha, \beta\ (\alpha,\beta \in {P})$
by 
\bea
&&~[h_{\bar{\epsilon}_j},\bar{\epsilon}_k]
=\langle \bar{\epsilon}_j,
\bar{\epsilon}_k \rangle,\qquad [h_{\bar{\epsilon}_j},h_{\bar{\epsilon}_k}]=0=
[\bar{\epsilon}_j,\bar{\epsilon}_k],\lb{HA1}
\ena
$h_\alpha=\sum_j n_j h_{\bep_j}$ for $\alpha=\sum_j n_j \bep_j$ and $h_0=0$.
Note that $\langle \bar{\epsilon}_j, \bar{\epsilon}_k \rangle=
\delta_{j,k}-\frac{1}{N}$ and $[h_{\alpha_j}, \alpha_k]=2\delta_{j,k}-\delta_{j,k+1}-\delta_{j,k-1}=A_{j k}$. 
Hence, we identify $h_{\alpha_j}=-h_{\bar{\epsilon}_j}+
h_{\bar{\epsilon}_{j+1}}$ with $h_j$ 
in the Drinfeld generators of $U_{q}(\slnh)$ (Section \ref{drinfelduq}).
 Noting $\sum_{j=1}^N h_{\bep_j}=0$, one can solve a set of equation $h_j=-h_{\bar{\epsilon}_j}+h_{\bar{\epsilon}_{j+1}}\ (1\leq j\leq N-1)$ for $h_{\bar{\epsilon}_j}$.
\bea 
&&h_{\bar{\epsilon}_j}=\sum_{k=1}^{j-1} h_k - \frac{1}{N}
\sum_{k=1}^{N-1} (N-k) h_k.\lb{hbepj}
\ena
From this and \eqref{def:D2}-\eqref{def:D4}, one can verify
the following commutation relations with 
the Drinfeld generators $c,\ d,\ a_{j,m}, 
x_{j,m}^\pm$ of $U_q(\widehat{\goth{sl}}_N)$.
\begin{eqnarray}
&&[h_{\bar{\epsilon}_i},a_{j_m}]=
[h_{\bar{\epsilon}_i},d]=
[h_{\bar{\epsilon}_i},c]=0,\\
&&~[h_{\bar{\epsilon}_i},x_{j,m}^\pm]=\pm(-\delta_{i,j}+
\delta_{i,j+1})x_{j,m}^\pm.
\end{eqnarray}

Now let us introduce another Heisenberg algebra $\C\{\H\}$ generated by
$P_{{\alpha}}$ and $ Q_{{\beta}}\ 
({\alpha}, {\beta} \in {P})$ satisfying the 
commutation relations
\begin{eqnarray}
&&[P_{\bar{\epsilon}_j}, Q_{\bar{\epsilon}_k}]=
\langle \bar{\epsilon}_j, \bar{\epsilon}_k \rangle, \qquad 
[P_{\bar{\epsilon}_j}, P_{\bar{\epsilon}_k}]=0=
[Q_{\bar{\epsilon}_j}, Q_{\bar{\epsilon}_k}],\lb{HA2}
\end{eqnarray}
where
$P_\alpha=\sum_j n_j P_{\bep_j}$ for $\alpha=\sum_j n_j \bep_j$ and $P_0=0$.
We also impose that $\C\{\H\}$ commutes with $U_q(\slnh)$.
\begin{eqnarray}
&&[P_{\bar{\epsilon}_j}, \al]=
[Q_{\bar{\epsilon}_j}, \al]=0, \lb{HA3}\\
&&[P_{\bar{\epsilon}_j}, U_q(\widehat{\goth{sl}}_N)]=
[Q_{\bar{\epsilon}_j}, U_q(\widehat{\goth{sl}}_N)]=0.\lb{HA4} 
\end{eqnarray}


\begin{dfn}
We define an extension $\C\{\hat{\H}\}$ of the Heisenberg algebra $\C\{\H \}$ 
by introducing new generators $\bbep_j\ (1\leq j\leq N-1)$
 and modifying the relations \eqref{HA2} to the following ones.
\bea
&&[P_{\bar{\epsilon}_j}, Q_{\bar{\epsilon}_k}]=
\langle \bar{\epsilon}_j, \bar{\epsilon}_k \rangle, \qquad 
[P_{\bar{\epsilon}_j}, P_{\bar{\epsilon}_k}]=0,
\lb{central1}
\\
&&
[Q_{\bar{\epsilon}_{j}}, Q_{\bar{\epsilon}_{k}}]=
\left(\frac{1}{r}-\frac{1}{r^*}\right){\rm sgn}(j-k)
{\rm log}~q,\lb{central2}\\
&&[Q_{\bar{\epsilon}_{j}},\bbep_k]=
\frac{1}{r}{\rm sgn}(j-k)\log~q,\lb{central3}\\
&&[\bbep_j,\bbep_k]=\frac{1}{r}{\rm sgn}(j-k)\log~q,\lb{central4}\\
&&[P_{\bar{\epsilon}_j}, \bbep_k]=0,\qquad \sum_{j=1}^N \bbep_j=0.\lb{central5}
\ena
We also impose the following commutation
relations.
\bea
&&
[\bbep_j, \al]=
[\bbep_j, U_{q}(\slnh)]=0.
\lb{central6}
\ena
\end{dfn}
If we set $\ba_j=-\bbep_j+\bbep_{j+1}$, we have 
\begin{prop}
\bea
&&[Q_{\al_{j}}, Q_{\al_{k}}]=
\left(\frac{1}{r}-\frac{1}{r^*}\right)\left(\delta_{j,k+1}-
\delta_{j,k-1}\right)
{\rm log}~q,\lb{qajqak}\\
&&[Q_{\bep_{j}}, Q_{\al_{k}}]=-
\left(\frac{1}{r}-\frac{1}{r^*}\right)\left(\delta_{j,k}+
\delta_{j,k+1}\right)
{\rm log}~q,\lb{qbejqak}\\
&&[Q_{\bar{\epsilon}_{j}},\bar{\alpha}_k]=
-\frac{1}{r}\left(\delta_{j,k}+\delta_{j,k+1}\right)
{\rm log}~q,\lb{qbejbak}\\
&&[Q_{\al_{j}},\bar{\alpha}_k]=
\frac{1}{r}\left(\delta_{j,k+1}-\delta_{j,k-1}\right)
{\rm log}~q,\lb{qajbak}\\
&&[\bar{\alpha}_j, \bar{\alpha}_k]=\frac{1}{r}\left(\delta_{j,k+1}-
\delta_{j,k-1}\right){\rm log}~q,
\lb{bajbak}\\
&&[\ba_j, P_{\bar{\epsilon}_k} ]=[\ba_j, U_q(\slnh)]=0.\lb{pbejbak}
\ena
\end{prop}

\subsubsection{Definition of $U_{q,p}(\slnh)$}
Now we are ready to define the currents 
$E_j(v), F_j(v), H_j^\pm(v)\ (1\leq j\leq N-1)$ and $K_j(v)\ (1\leq j\leq N)$.
\begin{eqnarray}
&&{E}_j(v)=e_j(z,p)e^{\bar{\alpha}_j}e^{-Q_{\alpha_j}}(q^{-j+N}z)
^{-\frac{P_{\alpha_j}-1}{r^*}},\label{def:total1c}\\
&&{F}_j(v)=f_j(z,p)e^{-\bar{\alpha}_j}
(q^{-j+N}z)^{\frac{P_{\alpha_j}-1}{r}}
(q^{-j+N}z)^{\frac{h_j}{r}},\label{def:total2c}\\
&&{H}_j^\pm
\left(v\right)=
\psi_j^\pm(z,p)
q^{\mp h_j}e^{-Q_{\alpha_j}}
(q^{-j+N\pm(r-\frac{c}{2})}z)^{(-\frac{1}{r^*}+\frac{1}{r})
(P_{\alpha_j}-1)+\frac{1}{r}h_j},\label{def:total3c}
\\
&&{K}_j(v)=k_j(z,p)e^{Q_{\bar{\epsilon}_j}}
z^{(\frac{1}{r^*}-\frac{1}{r})P_{\bar{\epsilon}_j}}
z^{-\frac{1}{r}h_{\bar{\epsilon}_j}+(\frac{1}{r^*}-\frac{1}{r})\frac{N-1}{2N}}.\label{def:total4c}
\end{eqnarray}
Here the currents $e_j(z,p), f_j(z,p), \psi_j^\pm(z,p)$ 
and $k_j(z,p)$ are the elliptic currents of $U_q(\slnh)$ 
given in (\ref{def:e})-(\ref{def:psi-}) and (\ref{def:k}), whereas 
$\ba_j,\ P_{\alpha},\ Q_{\beta}\ (\al,\ \beta \in P)$ are the elements in 
the Heisenberg algebra $\C\{\hat{\H}\}$.
From \eqref{ee}-\eqref{ef}, \eqref{kjkj}-\eqref{kifj} and 
 \eqref{central1}-\eqref{central6},
we can verify  the following relations.
\begin{prop}\lb{uqpslnh}
\begin{eqnarray}
&&E_i(v_1)E_j(v_2)=
\frac{[v_1-v_2+\frac{A_{ij}}{2}]^*}{[v_1-v_2-\frac{A_{ij}}{2}]^*}
E_j(v_2)E_i(v_1),
\label{EiEj}\\
&&\nn\\
&&F_i(v_1)F_j(v_2)=\frac{[v_1-v_2-\frac{A_{ij}}{2}]}
{[v_1-v_2+\frac{A_{ij}}{2}]}
F_j(v_2)F_i(v_1),
\label{FiFj}\\
&&\nn\\
&&[E_i(v_1),F_j(v_2)]=\frac{\delta_{i,j}}{q-q^{-1}}
\left(\delta(q^{-c}z_1/z_2)H_j^+\left(v_2+\frac{c}{4}\right)
-\delta(q^{c}z_1/z_2)H_j^-\left(v_2-\frac{c}{4}\right)
\right),
\label{EiFj}\\
&&\nn\\
&&H_j^\pm\left(v\mp \frac{1}{2}\left(r-\frac{c}{2}\right)\right)
=\kappa K_j\left(v+\frac{N-j}{2}\right)
K_{j+1}\left(v+\frac{N-j}{2}\right)^{-1},\label{HKK}\\
&&\nn\\
&&{K_j}(v_1){K_j}(v_2)={\rho^{}(v_1-v_2)}{K_j}(v_2){K_j}(v_1)
,\label{k2}
\\
&&{K_{j_1}}(v_1){K_{j_2}}(v_2)=
{\rho(v_1-v_2)}
\frac{[v_1-v_2-1]^*[v_1-v_2]}{[v_1-v_2]^*[v_1-v_2-1]}
{K_{j_2}}(v_2){K_{j_1}}(v_1)
\nonumber\\
&&~~~~~~~~~~~~~~~~\qquad\qquad (1\leqq j_1<j_2\leqq N),\label{k3}\\
&&{K_j}(v_1){E_j}(v_2)=
\frac{[v_1-v_2+\frac{j+r^*-N}{2}]^*}
{[v_1-v_2+\frac{j+r^*-N}{2}-1]^*}
{E_j}(v_2){K_j}(v_1),\label{k4}\\
&&{K_{j+1}}(v_1){E_j}(v_2)=
\frac{[v_1-v_2+\frac{j+r^*-N}{2}]^*}
{[v_1-v_2+\frac{j+r^*-N}{2}+1]^*}
{E_j}(v_2){K_{j+1}}(v_1),\label{k5}\\
&&{K_{j_1}}(v_1){E_{j_2}}(v_2)=
{E_{j_2}}(v_2){K_{j_1}}(v_1)~~\qquad (j_1\neq j_2,j_2+1),
\label{k6}\\
&&\nn\\
&&{K_j}(v_1){F_j}(v_2)=
\frac{[v_1-v_2+\frac{j+r-N}{2}-1]}
{[v_1-v_2+\frac{j+r-N}{2}]}
{F_j}(v_2){K_j}(v_1),
\label{k7}
\\
&&K_{j+1}(v_1)F_j(v_2)=
\frac{[v_1-v_2+\frac{j+r-N}{2}+1]}
{[v_1-v_2+\frac{j+r-N}{2}]}F_j(v_2)K_{j+1}(v_1)
,\label{k8}
\\
&&K_{j_1}(v_1)F_{j_2}(v_2)=
F_{j_2}(v_2)K_{j_1}(v_1)~~\qquad (j_1\neq j_2,j_2+1),\label{k9}\\
&&\nn\\
%
&&z_1^{-\frac{1}{r^*}}
\frac{(p^*q^2z_2/z_1;p^*)_\infty}
{(p^*q^{-2}z_2/z_1;p^*)_\infty
}\left\{
({z_2}/{z})^{\frac{1}{r^*}}
\frac{(p^*q^{-1}z/z_1;p^*)_\infty 
(p^*q^{-1}z/z_2;p^*)_\infty}
{(p^*qz/z_1;p^*)_\infty 
(p^*qz/z_2;p^*)_\infty}E_i(v_1)E_i(v_2)E_j(v)\right.\nonumber\\
&&
-\left.[2]_q\frac{(p^*q^{-1}z/z_1;p^*)_\infty 
(p^*q^{-1}z_2/z;p^*)_\infty}
{(p^*qz/z_1;p^*)_\infty 
(p^*qz_2/z;p^*)_\infty}
E_i(v_{1})E_j(v)E_i(v_{2})
\right.\nn
\\
&&+\left.
(z/z_1)^{\frac{1}{r^*}}\frac{(p^*q^{-1}z_1/z;p^*)_\infty 
(p^*q^{-1}z_2/z;p^*)_\infty}
{(p^*qz_1/z;p^*)_\infty 
(p^*qz_2/z;p^*)_\infty}
E_j(v)E_i(v_{1})E_i(v_{2})
\right\}+(z_1 \leftrightarrow z_2)=0,\nn\\
&&\label{e20}\\
&&z_1^{\frac{1}{r}}
\frac{(pq^{-2}z_2/z_1;p)_\infty}
{(pq^{2}z_2/z_1;p)_\infty
}\left\{
(z/z_2)^{\frac{1}{r}}
\frac{(pq z/z_1;p)_\infty 
(pq z/z_2;p)_\infty}
{(pq^{-1} z/z_1;p)_\infty 
(pq^{-1} z/z_2;p)_\infty}
F_i(v_1)F_i(v_2)F_j(v)\right.\nonumber\\
&&
-\left.[2]_q\frac{(pq z/z_1;p)_\infty 
(pq z_2/z;p)_\infty}
{(pq^{-1} z/z_1;p)_\infty 
(pq^{-1} z_2/z;p)_\infty}
F_i(v_{1})F_j(v)F_i(v_{2})
\right.\nn
\\
&&+\left.
(z_1/z)^{\frac{1}{r}}
\frac{(pq z_1/z;p)_\infty 
(pq z_2/z;p)_\infty}
{(pq^{-1} z_1/z;p)_\infty 
(pq^{-1} z_2/z;p)_\infty}
F_j(v)F_i(v_{1})
F_i(v_{2})
\right\} + (z_1 \leftrightarrow z_2)=0
\quad (|i-j|=1).\nonumber\\
&&~~~~~~~~~~~~~~~~~\label{e21}
\end{eqnarray}
Here the constant $\kappa$ and the function  $\rho(v)$ are given in
\eqref{psikk} and \eqref{rhofunc},
respectively. 
\end{prop}
The following relations among $H_j^\pm(v)$ are also useful.
\begin{prop}
\bea
&&H_j^+\left(v-r+\frac{c}{4}\right)=
H_j^-\left(v-\frac{c}{4}\right),\label{e8}\\
&&H_i^+(v_1)H_j^+(v_2)=
\frac{[v_1-v_2-\frac{A_{ij}}{2}]
[v_1-v_2+\frac{A_{ij}}{2}]^*
}{[v_1-v_2+\frac{A_{ij}}{2}]
[v_1-v_2-\frac{A_{ij}}{2}]^*
}H_j^+(v_2)H_i^+(v_1),\label{HiHj}
\\
&&\nn\\
&&H_i^+(v_1)E_j(v_2)=\frac{[v_1-v_2+\frac{A_{ij}}{2}+\frac{c}{4}]^*}
{[v_1-v_2-\frac{A_{ij}}{2}+\frac{c}{4}]^*}E_j(v_2)
H_i^+(v_1)
,\label{HiEj}
\\
&&\nn\\
&&H_i^+(v_1)F_j(v_2)=\frac{[v_1-v_2-\frac{A_{ij}}{2}-\frac{c}{4}]}
{[v_1-v_2+\frac{A_{ij}}{2}-\frac{c}{4}]}F_j(v_2)
H_i^+(v_1)
,\label{HiFj}
\ena
\end{prop}

\begin{dfn}{\bf ( Elliptic algebra $U_{q,p}(\slnh)$ )}
We define the elliptic algebra $U_{q,p}(\slnh)$ to be the 
associative algebra of the currents $E_j(v),\ F_j(v)\ (1\leq j\leq N-1)$ and 
$ K_j(v)\ (1\leq j\leq N)$ satisfying the relations 
\eqref{EiEj}-\eqref{e21}.
\end{dfn}

\begin{prop}
The construction of $E_j(v),\ F_j(v)$ and 
$ K_j(v)$ given in \eqref{def:total1c}-\eqref{def:total4c}
is a realization of
the elliptic algebra $U_{q,p}(\slnh)$
in terms of the Drinfeld generators of $U_q(\slnh)$,  
$h_{\bep_j},\  B^j_{m},\ x_{i,n}^\pm$ $(i=1,\cdots,N-1: m 
\in {{\Z}_{\neq 0}},\ n \in {\mathbb{Z}}),\ c,\ d$ and the Heisenberg algebra
$\C\{\hat{\H}\}$ generated by 
$P_{\bep_j},\ Q_{\bep_j}\ (1\leq j\leq N)$ and  $\ba_j\  (1\leq j\leq N-1)$.  
\end{prop}

\noindent
{\it Remark.}\ 
In Appendix A of \cite{JKOS2}, a realization of 
the elliptic algebra $U_{q,p}(\slnh)$ is given   
by using the Drinfeld currents of 
$U_q(\widehat{\goth{sl}}_N)$ and the Heisenberg algebra generated by 
$\{P_{j}, Q_{j}\}$ satisfying 
\bea
&&[P_{i}, Q_{j}]=-\frac{A_{ij}}{2},
\ena
which has no central extension.
The relation between $\{P_{j}, Q_{j}\}$ 
and $\{P_{\alpha_j}, Q_{\alpha_j}\}$ in $\C\{\H\}$ is
 $P_{\alpha_j}=P_{j}$ and
$Q_{\alpha_j}=-2Q_{j}$. The role of the central extension 
and the additional elements $\eta_j$ \eqref{central2}-\eqref{central4}
 is to suppress some extra 
$q$-fractional-power-factors in the relations in Proposition \ref{uqpslnh}. 
As for the problem realizing the 
$L$-operators satisfying the dynamical $RLL$-relation
 (Section \ref{loperator}), 
such $q$-factors can be absorbed into a choice of the gauge expressing 
the $R$-matrix. Conversely, in a gauge expressing the $R$-matrix 
components 
as $b(u,s)=q^{-\frac{1}{r}}
\frac{[s+1][s-1][u]}{[s]^2 [u+1]}$, 
$\bar{b}(u)=q^{\frac{1}{r}}\frac{[u]}{[u+1]}$ and the others remaining
 the same as in \eqref{rmatcomp}, we need neither the
central extension nor the addition of $\eta_j$.

\section{Half Currents}\lb{halfcurrents}
In order to construct a $L$-operator, 
we here introduce 
the half currents
$E_{l,j}^+(v),F_{j,l}^+(v)$ and $K_j^+(v)$ and investigate their 
commutation relations. We follow the idea of \cite{DF,EF,JKOS2}.

We often use the abbreviations
\begin{eqnarray}
P_{j,l}&=&
-P_{\bar{\epsilon}_{j}}+P_{\bar{\epsilon}_{l}}
=P_{\alpha_{j}}+P_{\alpha_{j+1}}+\cdots
+P_{\alpha_{l-1}},\\
h_{j,l}&=&-h_{\bar{\epsilon}_j}+h_{\bar{\epsilon}_l}=h_{j}+h_{j+1}+\cdots+
h_{l-1}~~
\end{eqnarray}
for $ j<l $. From the definition of $\C\{\hat{\H}\}$ and \eqref{def:total1c}-\eqref{def:total4c}, we have
\begin{eqnarray}
&&~[K_j(v),P_{k,l}]=(\delta_{j,k}-\delta_{j,l})
K_j(v)=[K_j(v),P_{k,l}+h_{k,l}],\lb{kjpjl}\\
&&~[E_j(v),P_{k,l}]=(\delta_{j,k}+
\delta_{j+1,l}-\delta_{j,l}-\delta_{j+1,k})E_j(v),\lb{ejpjl}\\
&&~[F_j(v),P_{j,l}+h_{j,l}]=(
\delta_{j,k}+
\delta_{j+1,l}-\delta_{j,l}-\delta_{j+1,k}
)F_j(v),\lb{fjpjlphjl}\\
&&
~[F_j(v),P_{k,l}]=0=
[E_j(v),P_{k,l}+h_{k,l}].\lb{kakan}
\end{eqnarray}

Now we define the half currents of $U_{q,p}(\slnh)$
as follows.
\begin{df}~~{\bf (Half currents)}~~
We define the half currents
$F_{j,l}^+(v), E_{l,j}^+(v), (1\leq j<l \leq N)$
and $K_j^+(v)\ (j=1,\cdots,N)$
by
\begin{eqnarray}
K_j^+(v)&=&K_j\left(v+\frac{r+1}{2}\right)~~~
(1\leq j \leq N),\lb{kjp}\\
F_{j,l}^+(v)&=&a_{j,l}\oint_{C(j,l)}
\prod_{m=j}^{l-1}\frac{dz_m}{2\pi i z_m}
F_{l-1}(v_{l-1})F_{l-2}(v_{l-2})
\cdots F_{j}(v_{j})\nonumber\\
&&\qquad\times
\frac{[v-v_{l-1}+P_{j,l}+h_{j,l}+\frac{l-N}{2}-1]
[1]}{[v-v_{l-1}+\frac{l-N}{2}]
[P_{j,l}+h_{j,l}-1]}\nonumber\\
&&\qquad\times
\prod_{m=j}^{l-2}
\frac{[v_{m+1}-v_{m}+P_{j,m+1}+h_{j,m+1}-\frac{1}{2}][1]}
{[v_{m+1}-v_{m}+\frac{1}{2}][P_{j,m+1}+h_{j,m+1}
]},\lb{fjl}\\
E_{l,j}^+(v)&=&a_{j,l}^*\oint_{C^*(j,l)}
\prod_{m=j}^{l-1}\frac{dz_m}{2\pi i z_m}
E_{j}(v_{j})E_{j+1}(v_{j+1})\cdots E_{l-1}(v_{l-1})
\nonumber\\
&&\qquad\times
\frac{[v-v_{l-1}-P_{j,l}+\frac{l-N}{2}+
\frac{c}{2}+1]^*[1]^*}
{[v-v_{l-1}+\frac{l-N}{2}+\frac{c}{2}]^*[P_{j,l}-1]^*
}
\nonumber\\
&&\qquad\times
\prod_{m=j}^{l-2}
\frac{[v_{m+1}-v_{m}-P_{j,m+1}+\frac{1}{2}]^*[1]^*}
{[v_{m+1}-v_{m}+\frac{1}{2}]^*[P_{j,m+1}-1]^*}.\lb{elj}
\end{eqnarray}  
Here the integration contour $C(j,l)$ and $C^*(j,l)$ are given by
\bea
C(j,l)&:& |pq^{l-N}z|<|z_{l-1}|<|q^{l-N}z|,\nn\\
&& |pqz_{k+1}|<|z_{k}|<|qz_{k+1}|,\\
C^*(j,l)&:& |p^*q^{l-N+c}z|<|z_{l-1}|<|q^{l-N+c}z|,\nn\\
&& |p^*qz_{k+1}|<|z_{k}|<|qz_{k+1}|,
\ena
where $k=j,j+1,..,l-2$.
The constants $a_{j,l}$ and $a_{j,l}^*$ are chosen to 
satisfy
\begin{eqnarray}
&&\frac{\kappa\ a_{j,l}a_{j,l}^* [1]}{q-q^{-1}}=1.
\end{eqnarray}
\end{df}
Then we can verify the following commutation relations 
\begin{thm}\lb{KEF}~~~The half currents $E_{l,j}^+(v), 
F_{j,l}^+(v)$ and
$K_j(v)\ (1\leq j<l\leq N)$ satisfy the following relations. 
\begin{eqnarray}
&&K_j^+(v_1)K_j^+(v_2)=
\rho(v)K_j^+(v_2)K_j^+(v_1)
\qquad (1\leq j \leq N),\lb{KjKj}\\
&&K_{j}^+(v_1)K_{l}^+(v_2)=
\rho(v)
\frac{[v-1]^*[v]}{[v]^*[v-1]}
K_{l}^+(v_2)K_{j}^+(v_1) \qquad
(1\leq j<l \leq N),\lb{kjkl}\\
&&\nn\\
&&K_l^+(v_1)^{-1}E_{l,j}^+(v_2)K_l^+(v_1)
=E_{l,j}^+(v_2)\frac{1}{\bar{b}_{}^*(v)}
-E_{l,j}^+(v_1)\frac{c_{}^*(v,P_{j,l})}{\bar{b}_{}^*(v)},\lb{klelj}\\
%
%
%
%
&&K_l^+(v_1)F_{j,l}^+(v_2)
K_l^+(v_1)^{-1}=
\frac{1}{\bar{b}_{}(v)}F_{j,l}^+(v_2)-
\frac{\bar{c}_{}(v,P_{j,l}+h_{j,l})}{
\bar{b}_{}(v)}F_{j,l}^+(v_1),\lb{klfjl}\\
%
%
%
%
&&\nn\\
&&\frac{[1-v]^*}{[v]^*}E^+_{l,j}(v_1)E^+_{l,j}(v_2)+\frac{[1+v]^*}{[v]^*}
E^+_{l,j}(v_2)E^+_{l,j}(v_1)\nn\\
&&\qquad\qquad=E^+_{l,j}(v_1)^2\frac{[1]^*[P_{j,l}-2+v]^*}{[P_{j,l}-2]^*[v]^*}+
E^+_{l,j}(v_2)^2\frac{[1]^*[P_{j,l}-2-v]^*}{[P_{j,l}-2]^*[v]^*},\lb{eljelj}\\
%
%
%
%
&&\nn\\
&&\frac{[1+v]}{[v]}F^+_{j,l}(v_1)F^+_{j,l}(v_2)+\frac{[1-v]}{[v]}
F^+_{j,l}(v_2)F^+_{j,l}(v_1)\nn\\
&&\qquad\qquad=F^+_{j,l}(v_1)^2\frac{[1][P_{j,l}+h_{j,l}-2-v]}{[P_{j,l}+h_{j,l}-2][v]}+
F^+_{j,l}(v_2)^2\frac{[1][P_{j,l}+h_{j,l}-2+v]}{[P_{j,l}+h_{j,l}-2][v]},
\lb{fjlfjl}\\
%
%
%
%
&&\nn\\
&&K_l^+(v_2)^{-1}E^+_{l,k}(v_1)K_l^+(v_2)E^+_{l,j}(v_2)\nn\\
&&\qquad\qquad=K_l^+(v_1)^{-1}E^+_{l,j}(v_2)K_l^+(v_1)E^+_{l,k}(v_1)\bar{R}^{*kj}_{kj}(v,P_{j,k})\nn\\
&&\qquad\qquad
+K_l^+(v_1)^{-1}E^+_{l,k}(v_2)K_l^+(v_1)E^+_{l,j}(v_1)\bar{R}^{*kj}_{jk}(v,P_{j,k}) \qquad (j\not=k), 
\lb{elkelj}\\
%
%
%
%
&&\nn\\
&&F^+_{k,l}(v_2)K_l^+(v_2)F^+_{j,l}(v_1)K_l^+(v_2)^{-1}\nn\\
&&\qquad\qquad=\bar{R}^{jk}_{jk}(v,P_{j,k}+h_{j,k})F^+_{j,l}(v_1)K_l^+(v_1)F^+_{k,l}(v_2)K_l^+(v_1)^{-1}\nn\\
&&\qquad\qquad
+\bar{R}^{kj}_{jk}(v,P_{j,k}+h_{j,k})F^+_{k,l}(v_1)K_l^+(v_1)F^+_{j,l}(v_2)K_l^+(v_1)^{-1}\qquad (j\not=k), 
\lb{fklfjl}\\
%
%
%
%
&&\nn\\
&&[E^+_{l,l-1}(v_1),F^+_{j,l}(v_2)]
=F^+_{j,l-1}(v_2)K^+_{l-1}(v_2)K^+_l(v_2)^{-1}
{[P_{l-1,l}-v-1]^*[1]^*\over[v]^*[P_{l-1,l}-1]^*} \nn \\
&&\qquad\qquad\qquad\qquad-F^+_{j,l-1}(v_1)K^+_l(v_1)^{-1}K^+_{l-1}(v_1)
{[P_{j,l}+h_{j,l}-v-1][1]\over [v_1-v_2][P_{j,l}+h_{j,l}-1]},\lb{elfjl}\\
&&\nn\\
&&[E^+_{l,j}(v_1),F^+_{l-1,l}(v_2)]
=K^+_{l-1}(v_2)K^+_l(v_2)^{-1}E^+_{l-1,j}(v_2)
\frac{[P_{j,l}-v-1]^*[1]^*}{[v]^*[P_{j,l}-1]^*} \nn \\
&&\qquad\qquad\qquad\qquad-K^+_l(v_1)^{-1}
K^+_{l-1}(v_1)E^+_{l-1,j}(v_1)\frac{[P_{l-1,l}+h_{l-1,l}-v-1][1]}
{[v][P_{l-1,l}+h_{l-1,l}-1]},\lb{eljfl}
\end{eqnarray}
where $v=v_1-v_2$.
\end{thm}
{\it Proof.}~~~The relations \eqref{KjKj} and \eqref{kjkl} are  direct 
consequences of \eqref{k2} and \eqref{k3}. 

We show the relation \eqref{klfjl}. The relations \eqref{klelj}
 can be proved in the same way. Setting $\pi_{l,j}=P_{j,l}+h_{j,l}$, we have 
from \eqref{k7}-\eqref{k8} and \eqref{kjpjl},
\bea
&&K_l^+(v_1)F_{j,l}^+(v_2)K_l^+(v_1)^{-1}\nn\\
&&=a_{j,l}\oint_{C(j,l)}
\prod_{k=j}^{l-1}\frac{dz'_k}{2\pi i z'_k}
F_{l-1}(v'_{l-1})F_{l-2}(v'_{l-2})
\cdots F_{j}(v'_{j})\nonumber\\
&&\times
\frac{[v_1-v'_{l-1}+\frac{l-N}{2}+1][v_2-v'_{l-1}+\pi_{l,j}+\frac{l-N}{2}-2]
[1]}{[v_1-v'_{l-1}+\frac{l-N}{2}][v_2-v'_{l-1}+\frac{l-N}{2}]
[\pi_{l,j}-2]}{\cal A}(v_{l-1},..,v_j;\pi_{l-1,j},..,\pi_{j+1,j}),\nn
\ena
where we set
\bea
{\cal A}(v'_{l-1},..,v'_j;\pi_{l-1,j},..,\pi_{j+1,j})=\prod_{k=j}^{l-2}
\frac{[v'_{k+1}-v'_{k}+\pi_{k+1,j}-\frac{1}{2}][1]}
{[v'_{k+1}-v'_{k}+\frac{1}{2}][\pi_{k+1,j}
]}.
\ena
Then the relation \eqref{klfjl} follows from the theta function identity 
\bea
&&\frac{[u_1+t][u_2+s]}{[u_1][u_2][s]}
=\frac{[u_1-u_2+t][u_2+s+t]}{[u_1-u_2][u_2][s+t]}+
\frac{[u_2-u_1+s][u_1+s+t][t]}{[u_2-u_1][u_1][s][s+t]} \lb{thetaid}
\ena
with the replacement $u_i=v_i-v_{l-1}'+\frac{l-N}{2} \ (i=1,2),\ s=\pi_{l,j}-2,\ t=1 $. 

Proofs of \eqref{eljelj}-\eqref{fjlfjl} and 
\eqref{elkelj}-\eqref{fklfjl} are lengthy. We put them in Appendix 
\ref{Proof}. 
 
Next let us consider the relation \eqref{elfjl}.
Integrating the delta function appearing from \eqref{EiFj}, we have 
\bea
&&(a_{j,l}a_{j,l}^*)^{-1}(q-q^{-1})[E^+_{l,l-1}(v_1),F^+_{j,l}(v_2)]\nn\\
&&=\Bigl\{\oint_{C^+_{l-1}}\du{z'_{l-1}}\cdots \oint\du{z'_{j}}\ 
H_{l-1}^+\left(v'_{l-1}+\frac{c}{4}\right)F^+_{l-1}(v'_{l-1})\cdots  F^+_{j}(v'_j)
\frac{[u_1-\pi_{l,l-1}+1]^*[1]^*}
{[u_1]^*[\pi_{l,l-i}-1]^*}\nn\\
&&-\oint_{C^-_{l-1}}\du{z'_{l-1}}\cdots \oint\du{z'_{j}}\ 
H_{l-1}^-\left(v'_{l-1}-\frac{c}{4}\right)F^+_{l-1}(v'_{l-1})\cdots  F^+_{j}(v'_j)
\frac{[u_1-\pi_{l,l-1}+1+c]^*[1]^*}
{[u_1+c]^*[\pi_{l,l-i}-1]^*}\Bigr\}\nn\\
&&\times {\cal A}(v'_{l-1},..,v'_j;\pi_{l-1,j},..,\pi_{j+1,j}).\nn
\ena
Here the contours $C^\pm_{l-1}$ are now 
\be
&& C^+_{l-1}\ {\rm encloses}\ z_1p^{*n}, z_2p^n\quad (n=1,2,...),\\
&&C^-_{l-1}\ {\rm encloses}\ z_1q^{2c}p^{*n}, z_2p^n\quad (n=1,2,...).
\en
Then in the second term, changing the variable $z'_{l-1}\to pz'_{l-1}$ 
and using the relation $H^-(v'+r-c/4)=H^+(v'+c/4)$,
we have the same integrand as the first term but the integration contour 
$C^-_{l-1}$ becomes 
\be
&&C^{-'}_{l-1}\ {\rm encloses}\ z_1p^{*n}, z_2p^n\quad (n=0,1,2,...).
\en
Therefore taking the residue at $z'_{l-1}=z_1, z_2$ and using the 
relation \eqref{HKK}, we get \eqref{elfjl}.

\begin{flushright}
{Q.E.D.}
\end{flushright}

\section{The $L$-operator of 
$U_{q,p}(\slnhbig)$ and Relation to $\Bqla(\slnhbig)$}\lb{loperator}
In this section, we construct a $L$-operator $\widehat{L}^+(u)$ 
 by using the half currents
and show that it satisfies the dynamical $RLL$-relation \eqref{DRLL2},
 which characterizes
the algebra ${\cal B}_{q,\lambda}(\slnh)$.
We then clarify the relation between the two elliptic algebras  
$U_{q,p}(\slnh)$ and ${\cal B}_{q,\lambda}(\slnh)$. 

\subsection{$L$-operator}
\begin{df}{\bf ($L$-operator)}~~By using the half currents, we define the $L$-operator
$\widehat{L}^+(v) \in {\rm End}({\mathbb{C}}^N)
\otimes U_{q,p}(\slnh)$
as follows.
\begin{eqnarray}
&&\widehat{L}^+(u)=
\left(\begin{array}{ccccc}
1&F_{1,2}^+(u)&F_{1,3}^+(u)&\cdots&F_{1,N}^+(u)\\
0&1&F_{2,3}^+(u)&\cdots&F_{2,N}^+(u)\\
\vdots&\ddots&\ddots&\ddots&\vdots\\
\vdots&&\ddots&1&F_{N-1,N}^+(u)\\
0&\cdots&\cdots&0&1
\end{array}\right)\left(
\begin{array}{cccc}
K^+_1(u)&0&\cdots&0\\
0&K^+_2(u)&&\vdots\\
\vdots&&\ddots&0\\
0&\cdots&0&K^+_{N}(u)
\end{array}
\right)\nn\\
&&\qquad\qquad\qquad\qquad\qquad\qquad\qquad\times
\left(
\begin{array}{ccccc}
1&0&\cdots&\cdots&0\\
E^+_{2,1}(u)&1&\ddots&&\vdots\\
E^+_{3,1}(u)&
E^+_{3,2}(u)&\ddots&\ddots&\vdots\\
\vdots&\vdots&\ddots&1&0\\
E^+_{N,1}(u)&E^+_{N,2}(u)
&\cdots&E^+_{N,N-1}(u)&1
\end{array}
\right).\lb{def:lhat}
\end{eqnarray}
Here $E_{l,j}^+(v), F_{j,l}^+(v)$ and
$K_j^+(v)$ are the half currents given in Section \ref{halfcurrents}.
\end{df}

Let $(\pi_{z}, V_z),\ V_z=V\otimes \C[z,z^{-1}]$ be the  evaluation
representation  
of $U_{q}(\slnh)$ based on the 
vector representation $V\cong \C^N$ (see Appendix \ref{Evaluation}). 
The image of the universal $R$-matrix $\cR(r,\{s_j\})$ of $\Bqla(\slnh)$ 
in the evaluation representation $(\pi_{V,z}\otimes \pi_{V,1})$ is 
given by the $R$-matrix $R^+(v,P)$ in \eqref{rmatfull}.
Then from a direct comparison with the relations 
of the half currents in Theorem \ref{KEF}, 
we conjecture the following property of the $L$-operator.

\noindent
\begin{conj}\lb{main}~~
The $L$-operator $\widehat{L}^+(v)$ satisfies
the following $RLL=LLR^*$ relation.
\begin{eqnarray}
R^{+(12)}(v_1-v_2,P+h)\widehat{L}^{+(1)}(v_1)
\widehat{L}^{+(2)}(v_2)=
\widehat{L}^{+(2)}(v_2)
\widehat{L}^{+(1)}(v_1)
R^{+*(12)}(v_1-v_2,P).\label{thm:RLL}
\end{eqnarray}
\end{conj}
In Appendix \ref{RLL}, we give a 
derivation of some of the relations among the half currents 
involved in \eqref{thm:RLL} and discuss their direct comparison with 
those in Theorem \ref{KEF}.
In Section 6.3, we give a proof of this statement in the case $c=1$.

\subsection{$U_{q,p}(\slnh)$ and $\Bqla(\slnh)$}

Based on the conjecture, we give a relation between $U_{q,p}(\slnh)$ and 
$\Bqla(\slnh)$. We argue that 
the $RLL$ relation \eqref{thm:RLL} is equivalent to
the dynamical $RLL$ relation of $\Bqla(\slnh)$. 
Hence we can regard the elliptic currents in $U_{q,p}(\slnh)$ as 
an elliptic analogue 
of the Drinfeld currents in $U_q(\slnh)$  providing 
a new realization of the elliptic quantum group $\Bqla(\slnh)$.

In order to show this, we consider the realization of $U_{q,p}(\slnh)$
 given in \eqref{def:total1c}-\eqref{def:total4c} and 
 modify the half currents in such a way that they have no 
 $Q_{\bar{\epsilon}_j},\ {\eta_j}\ (1\leq j \leq N)$ dependence. 
Let us define the modified half currents
$k_j^+(v,P)~(1\leq j \leq N)$
and  $e_{j,l}^+(v,P),\ f_{l,j}^+(v,P)~(1\leq j<l \leq N-1)$ as follows.
\begin{eqnarray}
k_j^+(v,P)&=&K_j^+(v)e^{-Q_{\bar{\epsilon}_j}},\\
e_{l,j}^+(v,P)&=&e^{Q_{\bar{\epsilon}_{l}}-{\bbep}_l}
E_{l,j}^+(v)e^{-Q_{\bar{\epsilon}_{j}}+{\bbep}_j},\\
f_{j,l}^+(v,P)&=&e^{-{\bbep}_j}F_{j,l}^+(v)e^{{\bbep}_l}.
\end{eqnarray}
Then it is easy to see from \eqref{def:total1c}-\eqref{def:total4c} and 
\eqref{kjp}-\eqref{elj}
that the modified half currents depend on neither 
$Q_{\bar{\epsilon}_j}$ nor $ {\bbep}_j$ and commute with 
$P_{\bar{\epsilon}_j}\ \forall j$. 
We hence regard them as the currents in $U_q(\slnh)$ with
parameters $P_{\bar{\epsilon}_j}$ and $r$. 

Now we define a modified $L$-operator ${L}^+(v,P)$ by
\begin{eqnarray}
&&{L}^+(u,P)=
\left(\begin{array}{ccccc}
1&f_{1,2}^+(u,P)&f_{1,3}^+(u,P)&\cdots&f_{1,N}^+(u,P)\\
0&  1           &f_{2,3}^+(u,P)&\cdots&f_{2,N}^+(u,P)\\
\vdots&\ddots    &\ddots        &\ddots&\vdots\\
\vdots&         &\ddots        &   1  &f_{N-1,N}^+(u,P)\\
0&\cdots        &\cdots        & 0    & 1
\end{array}\right)\nn\\
&&\times\left(
\begin{array}{cccc}
k^+_1(u,P)&0&\cdots&0\\
0&k^+_2(u,P)&&\vdots\\
\vdots&&\ddots&0\\
0&\cdots&0&k^+_{N}(u,P)
\end{array}
\right)
\left(
\begin{array}{ccccc}
1             &0&\cdots&\cdots&0\\
e^+_{2,1}(u,P)&1&\ddots&     &\vdots\\
e^+_{3,1}(u,P)&e^+_{3,2}(u,P)&\ddots&\ddots&\vdots\\
\vdots&\vdots&\ddots&1&0\\
e^+_{N,1}(u,P)&e^+_{N,2}(u,P)&\cdots&e^+_{N,N-1}(u,P)&1
\end{array}
\right).\nn\\
\lb{dL}
\end{eqnarray} 
Then the $L$-operator $\widehat{L}^+(v)$ and the modified one $L^+(v,P)$ 
are related by
\begin{eqnarray}
L^+(v,P)=
\widehat{L}^+(v)\left(\begin{array}{cccc}
e^{-Q_{\bar{\epsilon}_1}}&0&\cdots&0\\
0&e^{-Q_{\bar{\epsilon}_2}}&&\vdots\\
\vdots&&\ddots&0\\
0&\cdots&0&e^{-Q_{\bar{\epsilon}_N}}
\end{array}\right)=\widehat{L}^+(v)\ 
\exp\left\{\sum_{m=1}^N h^{(1)}_{{\epsilon}_m}Q_{\bar{\epsilon}_m}\right\}.
\lb{lhat}
\end{eqnarray}
Here $h^{(1)}_{{\epsilon}_j}=h_{{\epsilon}_j}\otimes 1,\ h_{\epsilon_m}\equiv
-E_{mm}$ (a $N \times N$ { matrix unit}). 
The reader should not confuse $h_{\epsilon_m}$ with $h_{\bep_m}$, but
note $h_j=-h_{\bep_j}+h_{\bep_{j+1}}=-h_{\ep_j}+h_{\ep_{j+1}}$ on $V$.

Substituting \eqref{lhat} into \eqref{thm:RLL} and 
noting the commutation relations
\begin{eqnarray}
&&P_{j,l}
\exp\left\{\sum_{m=1}^N h_{\epsilon_m}^{(k)}Q_{\bar{\epsilon}_m}\right\}
=\exp\left\{\sum_{m=1}^N h_{\epsilon_m}^{(k)}Q_{\bar{\epsilon}_m}\right\}
(P_{j,l}+h_{j,l}^{(k)})
\ena
and 
\bea
&&\left[ \sum_{m=1}^N( h^{(1)}_{{\epsilon}_m}+ h^{(2)}_{{\epsilon}_m})Q_{\bar{\epsilon}_m},\ R^{+*(12)}(v,P) \right]=0,
\ena 
or equivalently
\bea
&&\left[ Q_{\bep_j}+Q_{\bep_l},\ P_{j,l} \right]=0,
\ena
we can move each factor $\exp\left\{-\sum_{m=1}^N h^{(k)}_{{\epsilon}_m}
Q_{\bar{\epsilon}_m} \right\}\ (k=1,2)$ to the right end in the both sides. 
We then obtain the following statement.
\begin{cor}
The modified $L$-operator $L^+(v,P)$ satisfies
the dynamical $RLL$ relation
\begin{eqnarray}
&&R^{+(12)}(v,P+h)L^{+(1)}(v_1,P)L^{+(2)}(v_2,P+h^{(1)})
=
L^{+(2)}(v_2,P)L^{+(1)}(v_1,P+h^{(2)})
R^{+*(12)}(v,P),\nn\\
&&\label{thm:dRLL}
\end{eqnarray}
where $v=v_1-v_2$.
\end{cor}
Comparing this with \eqref{DRLL2}, 
we identify 
our $L^+(v,P)$ with $L^+(v,s)$ in \eqref{DRLL2} and 
$s_j$ with $P_{\al_j}$. Note the parametrization \eqref{sj}.
As a consequence of this result, we regard the elliptic currents 
$E_j(v),\ F_j(v)\ (1\leq j\leq N-1)$ and 
$K_j(v)\ (1\leq j\leq N)$ in $U_{q,p}(\slnh)$ 
as the Drinfeld currents of the elliptic 
quantum group $\Bqla(\slnh)$ up to tensoring with the Heisenberg algebra. 
Conversely,  this indicates that $U_{q,p}(\slnh)$ is  
an extension of the algebra $\Bqla(\slnh)$ 
by tensoring the Heisenberg algebra $\C\{\hat{\H}\}$ generated by 
$\{P_{\bep_j}, Q_{\bep_j}, \eta_j\}$. 
Namely, $U_{q,p}(\slnh)$ is obtained from 
${\Bqla}(\slnh)$, first by tensoring 
the half of the generators $e^{n Q_{\bep_j}}e^{m \eta_l} \ 
(1\leq j, l\leq N; n,m\in \Z)$, then
regarding $s_j=P_{\alpha_j}$ and imposing  
the commutation relations \eqref{central1}-\eqref{central6}. 
 Hence
\begin{eqnarray}
U_{q,p}(\slnh)=\Bqla(\widehat{\goth{sl}}_N)
\otimes_{{\mathbb{C}}\{P_{\bep_1},P_{\bep_2},..,P_{\bep_{N-1}}\}} 
{\mathbb{C}}\{\hat{\H}\}.
\end{eqnarray}


\section{Vertex Operators of $U_{q,p}(\slnhbig)$}
\lb{vertexoperators}

Tensoring the Heisenberg algebra breaks down 
the  coalgebra structure of $\Bqla(\slnh)$ \cite{JKOS2}. But 
we can still define $U_{q,p}(\slnh)$ counterparts of 
the intertwining operators of $\Bqla(\slnh)$. 
We call such operators the vertex operators of $U_{q,p}(\slnh)$.
In this section, we study
such vertex operators 
and compare them with those of
the $A_{N-1}^{(1)}$-type face model obtained in the papers
 \cite{AJMP,FKQ}.

\subsection{Intertwining relations}
Here we  derive $U_{q,p}(\slnh)$ counterparts of the dynamical 
intertwining relations \eqref{dintrelI}-\eqref{dintrelII}. 
In the next subsection, we use such relations to derive
a free field realization of the vertex operators. 

Let us first define an extension of the $U_q$ modules by
\be
&&\widehat{\cal F}=\bigoplus_{\mu_1,\cdots,\mu_N \in {\mathbb{Z}}}
{\cal F}\otimes e^{\mu_1 Q_{\bar{\epsilon}_1}
\cdots +\mu_N Q_{\bar{\epsilon}_N}}.
\en 
Let ${\Phi}_W(z,P)$ and ${\Psi}^*_W(z,P)$ be the 
type I and type II intertwining operators of $\Bqla(\slnh)$ \eqref{intI}
-\eqref{intII}.
We define type I and type II vertex operators $\widehat{\Phi}_W(v),\
 \widehat{\Psi}^*_W(v)$ of $U_{q,p}(\slnh)$
as the following extensions of the corresponding 
intertwining operators of $\Bqla(\slnh)$.
\begin{eqnarray}
&&\widehat{\Phi}_W(v)={\Phi}_W(q^cz,P)\qquad :
\widehat{\cal F} \longrightarrow \widehat{\cal F}'\otimes W_z,\\
&&\widehat{\Psi}^*_W(v)={\Psi}^*_W(z,P)
\exp\left\{ \sum_{j=1}^N h_{{\epsilon}_j} Q_{\bar{\epsilon}_j}\right\}\qquad :
W_z \otimes \widehat{\cal F} \longrightarrow \widehat{\cal F}'.
\end{eqnarray}
From the relations
\eqref{lhat} and \eqref{dintrelI}-\eqref{dintrelII},
the new operators $\widehat{\Phi}_W(v)$
and $\widehat{\Psi}_W^*(v)$
satisfy the following ``intertwining relations''. 
\begin{eqnarray}
\widehat{\Phi}_W^{(3)}(v_2)\widehat{L}_V^{+(1)}(v_1)&=&
R^{+(13)}_{VW}(v_1-v_2,P+h)\widehat{L}_V^{+(1)}(v_1)
\widehat{\Phi}_W^{(3)}(v_2),\label{Inter:Uqp1}\\
\widehat{L}_V^{+(1)}(v_1)\widehat{\Psi}_W^{*(2)}(v_2)&=&
\widehat{\Psi}_W^{*(2)}(v_2)\widehat{L}_V^{+(1)}(v_1)
R_{VW}^{+*(12)}(v_1-v_2,P-h^{(1)}-h^{(2)}).\label{Inter:Uqp2}
\end{eqnarray}

Now we restrict ourselves to the vector representation $V$ 
and investigate the relations \eqref{Inter:Uqp1}-\eqref{Inter:Uqp2}
 in detail.
We denote  a basis of $V$ by $\{\bfv_m \}_{m=1}^N$.
In this representation, the $R$-matrix $R^+_{VV}(v,P)$ is given by 
$R^+(v,P)$ in (\ref{rmatfull}) and the $L$-operator 
$\widehat{L}^+_V(v)$ by $\hL^+(v)$ in (\ref{def:lhat}).

We define the components of the vertex operators by
\begin{eqnarray}
\widehat{\Phi}_V\left(u-\frac{1}{2}\right)
=\sum_{m=1}^N \Phi_m(u)\otimes \bfv_m,~~
\widehat{\Psi}_V^*\left(u-\frac{c+1}{2}\right)(\bfv_m \otimes \cdot)=
\Psi^*_m(u),
\end{eqnarray}
and the matrix elements of
the $L$-operator $\widehat{L}^+(u)$ by
\begin{eqnarray}
\widehat{L}^+(u)\bfv_{j}=
\sum_{1\leq m \leq N}
\bfv_m L^+(u)_{mj}.
\end{eqnarray}

Using these components, the equation (\ref{Inter:Uqp1})
is read as follows.
\begin{eqnarray}
&&\Phi_m\left(v_2\right)
L^+_{mj}(v_1)
=\rho^+(v_1-v_2+1/2)L^+_{mj}(v_1)\Phi_m\left(v_2\right)
,\label{IntTypeI1}\\
&&\rho^+(v_1-v_2+1/2)^{-1}
\Phi_{m}\left(v_2\right)L^+_{lj}(v_1)\nn\\
&&\quad =b(v_1-v_2+1/2,P_{l,m}+h_{l,m})
L^+_{lj}(v_1)\Phi_{m}\left(v_2\right)\nn\\
&&\qquad\qquad\qquad\qquad\qquad+
c(v_1-v_2+1/2,P_{l,m}+h_{l,m})
L^+_{mj}(v_1)\Phi_{l}\left(v_2\right),
\label{IntTypeI2}\\
&&\rho^+(v_1-v_2+1/2)^{-1}
\Phi_{l}\left(v_2\right)L^+_{mj}(v_1)\nn\\
&&\quad =\bar{b}(v_1-v_2+1/2)
L^+_{mj}(v_1)\Phi_{l}\left(v_2\right)+
\bar{c}(v_1-v_2+1/2,P_{l,m}+h_{l,m})
L^+_{lj}(v_1)\Phi_{m}\left(v_2\right),\nn\\
\label{IntTypeI3}
\end{eqnarray}
for $1\leq l<m \leq N$ and $1\leq j \leq N$.
For the type II, we have the following set of the
equations arising from 
the equation (\ref{Inter:Uqp2}) 
\begin{eqnarray}
&&
L^+_{jm}(v_1)\Psi^*_m(v_2)=
\rho^{+*}(v_1-v_2+1)\Psi^*_m(v_2)L^+_{jm}(v_1),
\label{IntTypeII1}
\\
&&\rho^{+*}(v_1-v_2+1)^{-1}
L^+_{jl}(v_1)\Psi_{m}^*(v_2)\nn\\
&&\quad =
\Psi_{m}^*(v_2)L^+_{jl}(v_1)
b^*(v_1-v_2+1,P_{l,m})+
\Psi_{l}^*(v_2)L^+_{jm}(v_1)
\bar{c}^*(v_1-v_2+1,P_{l,m}),\label{IntTypeII2}
\\
&&\rho^{+*}(v_1-v_2+1)^{-1}
L^+_{jm}(v_1)\Psi_{l}^*(v_2)\nonumber\\
&&\quad =
\Psi_{l}^*(v_2)L^+_{jm}(v_1)
\bar{b}^*(v_1-v_2+1)+
\Psi_{m}^*(v_2)L^+_{jl}(v_1)
{c}^*(v_1-v_2+1,P_{l,m}),\label{IntTypeII3}
\end{eqnarray}
for $1\leq l<m \leq N$ and $1\leq j \leq N$.

Let us investigate equations 
(\ref{IntTypeI1})-(\ref{IntTypeI3}) 
in detail. From the component
$j=m=N$ of equation
(\ref{IntTypeI1}), we have
\begin{eqnarray}
\Phi_N(v_2)K_N^+(v_1)=\rho^+\left(v_1-v_2+\frac{1}{2}
\right)K_N^+(v_1)
\Phi_N(v_2).\label{VI1}
\end{eqnarray}
Setting  $1\leq j<m=N$ in (\ref{IntTypeI1}), we have
\begin{eqnarray}
\Phi_N(v_2)E^+_{N,j}(v_1)=E^+_{N,j}(v_1)\Phi_N(v_2). \label{VI2}
\end{eqnarray}
The following relations turn out to be  sufficient conditions 
for (\ref{VI2}) to hold.
\begin{eqnarray}
\Phi_N(v_2)E_j(v_1)=E_j(v_1)\Phi_N(v_2)~\qquad (1\leq j \leq N-1).
\label{VI3}
\end{eqnarray}

Next let us consider the component $l<m=j=N$ 
of equation (\ref{IntTypeI2}).
We set 
\bea
&&\rho^+(v)=-\frac{[v+1]}{\varphi(v)},
\label{rho}
\\
&&\varphi(v)
=(q^r z)^{\frac{N+1}{rN}}[v-1]\frac{\{pz\}\{pq^{2N}z\}\{q^{2N+2}/z\}
\{q^{-2}/z\}}{\{pq^{2N+2}z\}\{pq^{-2}z\}\{1/z\}
\{q^{2N}/z\}}. 
\ena
Then \eqref{IntTypeI2} with $m=j=N$ can be written as
\bea
&&\varphi(v_1-v_2+1/2)\ 
\Phi_{N}\left(v_2\right)F^+_{l,N}(v_1)K^+_N(v_1)\nn\\
&&\qquad =\frac{[P_{l,N}+h_{l,N}-1][P_{l,N}+h_{l,N}+1][v_1-v_2+1/2]}
{[P_{l,N}+h_{l,N}]^2}
F^+_{l,N}(v_1)K_N^+(v_1)\Phi_{m}\left(v_2\right)\nn\\
&&\qquad\qquad +
\frac{[P_{l,N}+h_{l,N}+v_1-v_2+1/2][1]}{[P_{l,N}+h_{l,N}]}
K^+_N(v_1)\Phi_{l}\left(v_2\right).\lb{IntTypeI2b}
\ena
In order to solve \eqref{IntTypeI2b}, 
let us assume that the operator product $K_N^+(v_1)\Phi_N(v_2)$
does not have a pole at $v_1-v_2+3/2+r=0$.
Later we will check  
that, for $c=1$, this assumption is satisfied in 
a free field realization.
Then from relations (\ref{VI1}) and (\ref{rho}), 
we conclude that the product $\Phi_N(v_2)K_N^+(v_1)$ in the LHS of 
\eqref{IntTypeI2b}
 has zero at $v_1-v_2+3/2+r=0$.
Therefore, setting  $v_1-v_2+3/2+r=0$ in (\ref{IntTypeI2b}), we have
\begin{eqnarray}
\Phi_l(v_2)&=&
K^+_N\left(v_1\right)^{-1}
\frac{[P_{j,N}+h_{j,N}+1]}{[P_{j,N}+h_{j,N}]}
F^+_{l,N}\left(v_1\right)K^+_N(v_1)\Phi_N(v_2)\nonumber\\
&=&
F^+_{l,N}\left(v_2-{1}/{2}-r\right)
\Phi_N(v_2) ~~\quad
(1 \leq l \leq N-1).\label{VI4}
\end{eqnarray}
Note that the shift of $v$ by $r$ in $F^+_{l,N}(v)$ yields a change of 
contour (see \eqref{Type-Icont}).
Substituting (\ref{VI1}) and (\ref{VI4}) into
(\ref{IntTypeI2}) for $l<m=j=N$, and using
Riemann's theta identity,
we find that \eqref{VI4} and the following relations are 
sufficient conditions 
for (\ref{IntTypeI2}) with $l<m=j=N$.
\begin{eqnarray}
&&F_{N-1}(v_1)\Phi_N(v_2)=
\frac{[v_1-v_2+\frac{1}{2}]}
{[v_1-v_2-\frac{1}{2}]}
\Phi_N(v_2)F_{N-1}(v_1),\label{VI5}\\
&&F_l(v_1)\Phi_N(v_2)=\Phi_N(v_2)F_l(v_1)~\qquad (1\leq l \leq N-2)
\label{VI6},\\
&&[\Phi_N(v),P_{j,k}+h_{j,k}]=
-\delta_{k,N}\Phi_N(v)~\qquad (j<k).
\label{VI7}
\end{eqnarray}
In the next section, we construct a free field realization of the type I
vertex operators using relations \eqref{VI1}, \eqref{VI3} 
and \eqref{VI4}-\eqref{VI7} for $c=1$. We then check that the resulting
vertex operators satisfy the remaining relations in \eqref{IntTypeI2} and
\eqref{IntTypeI3}.

Similarly, from the $j=m=N$ component 
of (\ref{IntTypeII1}), we have for the
type-II vertex operator
\begin{eqnarray}
K_N^+(v_1)\Psi_N^*(v_2)=
\rho^{+*}\left(v_1-v_2+1\right)
\Psi_N^*(v_2)K_N^+(v_1) \label{VII1}
\end{eqnarray}
and from the 
$1\leq j<m=N$ component of (\ref{IntTypeII1}),
\begin{eqnarray}
F^+_{j,N}(v_1)\Psi_N^*(v_2)=
\Psi_N^*(v_2)F^+_{j,N}(v_1)~~\qquad (1\leq j \leq N-1).
\label{VII2}
\end{eqnarray}
We find the following as  sufficient conditions for
(\ref{VII2}).
\begin{eqnarray}
F_{j}(v_1)\Psi_N^*(v_2)=
\Psi_N^*(v_2)F_{j}(v_1)~~\qquad (1\leq j \leq N-1).
\label{VII3}
\end{eqnarray}
To solve equation (\ref{IntTypeII2}) with 
$l<j=m=N$, we assume that 
the product 
$\Psi_N^*(v_2)
K_{N}^+(v_1)$ has no pole at $v_1-v_2+2+r^*=0$.  
Then the product
$K_N^+(v_1)\Psi_N^*(v_2)$ in the LHS has a zero 
at $v_1-v_2+2+r^*=0$ for the same reason as the type I case. 
Therefore, from \eqref{IntTypeII2} with 
$l<j=m=N$ and setting $v_1-v_2+2+r^*=0$, we have
\begin{eqnarray}
\Psi_l^*(v)=\Psi_N^*(v)E_{N,l}^+\left(v-\frac{c+1}{2}-r^*\right)
~~\qquad (1\leq l \leq N-1).
\label{VII4}
\end{eqnarray}
Then \eqref{VII4} and the following relations turns out to be the 
sufficient conditions for \eqref{IntTypeII2} and \eqref{IntTypeII3}.
\begin{eqnarray}
&&E_{N-1}(v_1)\Psi_N^*(v_2)=
\frac{[v_1-v_2-\frac{1}{2}]^*}
{[v_1-v_2+\frac{1}{2}]^*}
\Psi_N^*(v_2)E_{N-1}(v_1)\label{VII5},\\
&&E_j(v_1)\Psi_N^*(v_2)=\Psi_N^*(v_2)E_j(v_1)~\qquad
(1\leq j \leq N-2),\label{VII6}\\
&& [\Psi_N^*(v),P_{j,k}]=\delta_{k,N}\Psi_N^*(v)~\qquad 
(j<k).
\label{VII7}
\end{eqnarray}

\subsection{Free field realizations}
Now we construct a free
field realization of the vertex operators 
fixing the representation level 
$c=1$. 
For this purpose, we first consider the simple root operator 
$\al_j$ introduced in Section \ref{heisenberg}. We make the following 
 standard central extension. 
\bea
[\al_j, \al_k]=i\pi  A_{jk}. \lb{cext}
\ena
Setting $\hat{\al}_j=\al_j+\ba_j$ where $\ba_j$ is an element of 
the Heisenberg algebra $\C\{\hat{H}\}$, we have 
\bea
&&[\hat{\alpha}_j, \hat{\alpha}_k]=i\pi  A_{jk}
+\frac{1}{r}\left(\delta_{j,k+1}-
\delta_{j,k-1}\right){\rm log}~q,
\lb{bajbak}\\
&&[h_{\bep_j}, \hat{\al}_k]
=-\delta_{j,k}+\delta_{j,k+1},\\
&&[Q_{\bar{\epsilon}_{j}},\hat{\alpha}_k]=
-\frac{1}{r}\left(\delta_{j,k}+\delta_{j,k+1}\right)
{\rm log}~q,\lb{qbejbak}\\
&&[Q_{\al_{j}},\hat{\alpha}_k]=
\frac{1}{r}\left(\delta_{j,k+1}-\delta_{j,k-1}\right)
{\rm log}~q,\lb{qajbak}\\
&&[\hat{\al}_j, P_{\bar{\epsilon}_k} ]=0.\lb{pbejbak}
\ena
Then the following statement holds.
\begin{prop}~~
The currents $E_j(v)$ and $F_j(v)$ given by
\begin{eqnarray}
E_j(v)&=&:\exp\left(
-\sum_{m \neq 0}\frac{[r m]_q}{m[r^* m]_q}(-B_m^j+B_m^{j+1})
(q^{N-j}z)^{-m}\right):
e^{\hat{\alpha}_j}z^{h_j}e^{-Q_{\alpha_j}}(q^{-j+N}z)
^{-\frac{P_{\alpha_j}-1}{r^*}},\nonumber\\
\label{free:E}\\
F_j(v)&=&:\exp\left(
\sum_{m\neq 0}\frac{1}{m}(-B_m^j+B_m^{j+1})
(q^{N-j}z)^{-m}
\right):
e^{-\hat{\alpha}_j}z^{-h_j}(q^{-j+N}z)^{\frac{P_{\alpha_j}-1}{r}
+\frac{h_j}{r}},\label{free:F}
\end{eqnarray}
together with $H_j^\pm(v), K_j(v)$ given in
(\ref{def:total3c})-(\ref{def:total4c}) satisfy
the commutation relations in Proposition
\ref{uqpslnh} for level $c=1$. Hence they give
a free field realization of the level one 
elliptic algebra $U_{q,p}(\slnh)$.
\end{prop}

Now substituting 
the free field realization of
$E_j(v),\ F_j(v),\ K_j(v)$ into \eqref{kjp}-
\eqref{elj}, we
obtain a realization of 
the half currents 
$E^+_j(v),\ F^+_j(v),\ K^+_j(v)$ as well as 
the $L$-operator $\hL^+(v)$ satisfying the $RLL$-relation 
\eqref{thm:RLL} for $c=1$.
Using such a $L$-operator in the ``intertwining relations'',
 \eqref{VI1}-\eqref{VI7} for type I and \eqref{VII1}-\eqref{VII7} for the 
 type II, one can solve them for the vertex operators.
The results are stated as follows.    
\begin{thm}\lb{vertexop}
The highest components of the type I and the type II vertex operators 
$\Phi_N(v)$, $\Psi_N^*(v)$ are realized in terms of a free field by
\begin{eqnarray}
\Phi_N(v)&=&
:\exp\left(
-\sum_{m \neq 0}\frac{1}{m}B_m^N z^{-m}
\right):
e^{\bar{\Lambda}_{N-1}} 
z^{(1-\frac{1}{r})h_{\bar{\epsilon}_N}}
z^{-\frac{1}{r}P_{\bar{\epsilon}_N}}
z^{(1-\frac{1}{r})\frac{N-1}{2N}},\label{free:Phi}\\
\Psi_N^*(v)&=&:\exp\left(
\sum_{m \neq 0}\frac{[r m]_q}{m[r^* m]_q}B_m^N z^{-m}
\right):
e^{-\bar{\Lambda}_{N-1}}z^{-h_{\bar{\epsilon}_N}}
e^{Q_{\bar{\epsilon}_N}}
z^{\frac{1}{r^*}P_{\bar{\epsilon}_N}}
z^{(1+\frac{1}{r^*})\frac{N-1}{2N}}
q^{(\frac{1}{r}-\frac{1}{2})\frac{N-1}{N}}
,\nonumber\\
\label{free:Psi}
where
\end{eqnarray}
\begin{eqnarray}
\bar{\Lambda}_{N-1}&=&\frac{1}{N}
(\hat{\alpha}_1+2\hat{\alpha}_2+
\cdots +(N-1)\hat{\alpha}_{N-1}).
\end{eqnarray}
For the other components of the type I vertex operator $\Phi_j(v)\ (j=1,\cdots,N)$, we obtain from (\ref{VI4})
\begin{eqnarray}
\Phi_j(v)&=& a_{j,N}
\oint_{C} \prod_{m=j}^{N-1}
\frac{dz_m}{2\pi i z_m}
\Phi_N(v)F_{N-1}(v_{N-1})
\cdots F_{j}(v_{j})
\nonumber\\
&&\times
\prod_{m=j}^{N-1}
\frac{[v_{m+1}-v_{m}+
P_{j,m+1}+h_{j,m+1}
-\frac{1}{2}][1]}{[v_{m+1}-v_{m}+\frac{1}{2}]
[P_{j,m+1}+h_{j,m+1}]}\nn\\
&=&
a_{j,N}
\oint_{C} \prod_{m=j}^{N-1}
\frac{dz_m}{2\pi i z_m}
F_{j}(v_{j})\cdots
F_{N-1}(v_{N-1})
\Phi_N(v)
\nonumber\\
&&\times
\prod_{m=j}^{N-1}
\frac{[v_{m+1}-v_{m}+
P_{j,m+1}+h_{j,m+1}
-\frac{1}{2}][1]}{[v_{m+1}-v_{m}-\frac{1}{2}]
[P_{j,m+1}+h_{j,m+1}]},\label{Type-I}
\end{eqnarray}
where $v=v_N$ and  the integration contour $C$ 
is specified by the condition 
\begin{eqnarray}
&&|q^{-1}z|<|z_{N-1}|<|p^{-1}q^{-1}z|,\lb{Type-Icont}\\
&&
|pqz_{m+1}|<|z_{m}|<|q z_{m+1}|~~\qquad 
(j \leq m \leq N-2).\nn
\end{eqnarray}

For the type II vertex $\Psi_j^*(v)\ (j=1,\cdots,N)$, we obtain
from (\ref{VII4})
\begin{eqnarray}
\Psi_j^*(v)&=&
a_{j,N}^* \oint_{C^*} 
\prod_{m=j}^{N-1}
\frac{dz_m}{2\pi i z_m}
E_{j}(v_{j})\cdots 
E_{N-1}(v_{N-1})
\Psi_N^*(v_N)
\nonumber\\
&&\times
\prod_{m=j}^{N-1}
\frac{[v_{m+1}-v_{m}-
P_{j,m+1}
+\frac{1}{2}]^*[1]^*}
{[v_{m+1}-v_{m}+\frac{1}{2}]^*
[P_{j,m+1}-1]^*}\nn\\
&=&
a_{j,N}^* \oint_{C^*} 
\prod_{m=j}^{N-1}
\frac{dz_m}{2\pi i z_m}
\Psi_N^*(v_N)
E_{N-1}(v_{N-1})\cdots
E_{j}(v_{j})
\nonumber\\
&&\times
\prod_{m=j}^{N-1}
\frac{[v_{m+1}-v_{m}-
P_{j,m+1}
+\frac{1}{2}]^*[1]^*}
{[v_{m+1}-v_{m}-\frac{1}{2}]^*
[P_{j,m+1}-1]^*}.\label{Type-II}
\end{eqnarray}
The integration contour $C^*$ is specified as follows.
\begin{eqnarray}
|p^*q^{-1} z_{m+1}|,|q^{-1}z_{m+1}|
<|z_m|<|q z_{m+1}|,|p^{*-1}qz_{m+1}|~~\qquad(j \leq m \leq N-1).
\label{cont:3}
\end{eqnarray}
Here the integration variable $z_{m}\ (j\leq m\leq N-1)$ 
should encircle the poles 
$p^*q^{-1}z_{m+1}, q^{-1}z_{m+1}$ but not the poles 
 $p^{*-1} q z_{m+1}, qz_{m+1}$.
\end{thm}

In addition, 
we have the following commutation relations.
\begin{prop}~~~The highest components
$\Phi_N(v)$ and $\Psi_N^*(v)$
satisfy 
\begin{eqnarray}
&&[\Phi_N(v),P_{j_1,j_2}]=[\Psi_N^*(v),P_{j_1,j_2}+h_{j_1,j_2}]=0,\\
&&\Phi_N(v_1)\Psi_N^*(v_2)=\chi(v_1-v_2)
\Psi_N^*(v_2)\Phi_N(v_1), \label{VI-II}\\
&&\chi(v)=\frac{\Theta_{q^{2N}}(qz)}
{\Theta_{q^{2N}}(q/z)}.\lb{chifunc}
\end{eqnarray}
\end{prop}

\noindent
{\it Remark}~~~The free field realizations of the vertex operators in Theorem
\ref{vertexop} 
are essentially the same as those of the $A_{N-1}^{(1)}$-type face model 
obtained in \cite{AJMP,FKQ}. 
There are two differences between ours and those in \cite{AJMP,FKQ};
the choice of the gauge expressing the $R$-matrices and the zero-mode 
operators. Due to our gauge, we have the extra factors $\prod_{m=j}^{N-1}
\frac{[1]}{[P_{j,m+1}+h_{j,m+1}]}$  and 
$\prod_{m=j}^{N-1}\frac{[1]^*}{[P_{j,m+1}-1]^*}$ in the type I and the type II vertex operators, respectively. 
As for the zero-modes, the correspondence between ours $P_{\bar{\epsilon}_j},
Q_{\bar{\epsilon}_j},
h_j, \hat{\alpha}_j$ and those in 
\cite{AJMP, FKQ},  $P_{\alpha_j},
P_{\omega_N},
Q_{\alpha_j},
Q_{\omega_N}$ is given by 
\begin{eqnarray}
\left(\begin{array}{c}
\sqrt{\frac{r^*}{r}}P_{\alpha_j}
\\
\sqrt{\frac{r^*}{r}}P_{\omega_N}
\\
\sqrt{\frac{r}{r^*}}P_{\alpha_j}
\\
\sqrt{\frac{r}{r^*}}P_{\omega_N}
\end{array}\right)\leftrightarrow
\left(\begin{array}{c}
\frac{r^*}{r}h_j-\frac{1}{r}P_{\alpha_j}
\\
\frac{r^*}{r}h_{\bar{\epsilon}_N}
-\frac{1}{r}
P_{\bar{\epsilon}_N}\\
h_j-\frac{1}{r^*}P_{\alpha_j}\\
h_{\bar{\epsilon}_N}-\frac{1}{r^*}P_{\bar{\epsilon}_N}
\end{array}
\right)
,~~~~
\left(\begin{array}{c}
i\sqrt{\frac{r^*}{r}}Q_{\alpha_j}
\\
i\sqrt{\frac{r^*}{r}}Q_{\omega_N}
\\
i\sqrt{\frac{r}{r^*}}Q_{\alpha_j}
\\
i\sqrt{\frac{r}{r^*}}Q_{\omega_N}
\end{array}\right)
\leftrightarrow
\left(\begin{array}{c}
\hat{\alpha}_j\\
\bar{\Lambda}_{N-1}\\
\hat{\alpha}_j-Q_{\alpha_j}\\
\bar{\Lambda}_{N-1}-Q_{\bar{\epsilon}_N}
\end{array}\right).
\end{eqnarray}
One should note that 
to define the currents $K_j(v)$ by factoring the operators 
$H_j^\pm(v)$, the use of two sets of the Heisenberg operators $\{P_{\bar{\epsilon}_j}, Q_{\bar{\epsilon}_j}\}$ and $\{
h_j, \hat{\alpha}_j\}$ are essential.

\subsection{Commutation relations}
We next investigate commutation relations 
of the vertex operators and show that 
our realization satisfies the full intertwining relations for $c=1$.

\begin{thm}\label{thm:Vertexcom}~~~~
The free field realizations of
the type-I vertex operator
$\Phi_\mu(v)$ (\ref{Type-I}) and
the type-II vertex operator $\Psi_{\mu}^*(v)$
(\ref{Type-II})
satisfy the following commutation relations. 
\begin{eqnarray}
&&{\Phi}^{}_{j_2}(v_2)
{\Phi}^{}_{j_1}(v_1)=
\sum_{j_1',j_2'=1}^N
{R}_{j_1j_2}^{j_1'j_2'}(v_1-v_2,P+h)\ 
{\Phi}^{}_{j_1'}(v_1)
{\Phi}^{}_{j_2'}(v_2)
\label{Com:Type-I},\\
&&{\Psi}^{*}_{j_1}(v_1)
{\Psi}^{*}_{j_2}(v_2)
=\sum_{j_1',j_2'=1}^N
{\Psi}^{*}_{j_2'}(v_2)
{\Psi}^{*}_{j_1'}(v_1)\ 
{R}^{*j_1j_2}_{j_1'j_2'}(v_1-v_2,P)
\label{Com:Type-II},\\
&&{\Phi}^{}_{j}(v_1){\Psi}^{*}_{k}(v_2)
={\chi}(v_1-v_2)\ {\Psi}^{*}_{k}(v_2)
{\Phi}^{}_{j}(v_1).\label{Com:Type-I,II}
\end{eqnarray}
Here we set
\begin{eqnarray}
{R}(v,P+h)={\mu(v)}
\bar{R}(v,P+h),~~\quad
{R}^*(v,P)
={\mu^*(v)}
\bar{R}^{*}(v,P),
\end{eqnarray}
with
\begin{eqnarray}
&&\mu(v)=z^{({\frac{1}{r}-1)}\frac{N-1}{N}}
\frac{\{pq^{2N-2}z\}
\{q^2z\}
\{p/z\}
\{q^{2N}/z\}
}{
\{pz\}
\{q^{2N}z\}
\{pq^{2N-2}/z\}
\{q^2/z\}
},\\
\end{eqnarray}
and $\mu^*(v)=\mu(v)|_{r\to r^*}$. 
\end{thm}
{\it Proof.}~~~Using the formulae \eqref{VI4}-\eqref{VI7} and
 \eqref{VII4}-\eqref{VII7}, 
the commutation relations
(\ref{Com:Type-I}) and (\ref{Com:Type-II})
are reduced to the relations among the half currents 
 (\ref{eljelj2}), (\ref{fjlfjl2}), 
(\ref{elkelj}) and (\ref{fklfjl}).
Then the  proofs of the latter relations are given in Appendix \ref{Proof}.

Let us consider the relation (\ref{Com:Type-I,II}).
The case $j=N$ or $k=N$ 
is a direct consequence of \eqref{VI2}, \eqref{VII2} and  
\eqref{VI-II}.
The simplest non-trivial case is $j=k=N-1$.
From \eqref{VI4}, \eqref{VII4} and \eqref{EiFj}, 
we have the following equation after integrating
the delta functions.
\begin{eqnarray}
&&\Phi_{N-1}(v_1)\Psi_{N-1}^*(v_2)-
{\chi}\left(v_1-v_2\right)
\Psi_{N-1}^*(v_2)\Phi_{N-1}(v_1)\nonumber\\
&=&-\frac{a_{N-1,N}a_{N-1,N}^*}{q-q^{-1}}
\Phi_{N}(v_1)\Psi_{N}^*(v_2)\nn
\\
&\times&
\left(
\oint_{{C}_+} \frac{dz'}
{2\pi i z'} 
H_{N-1}^+\left(v'+\frac{1}{4}
\right)
\frac{[v_2-v'-P_{N-1,N}]^*[1]^*[v_1-v'+P_{N-1,N}+
h_{N-1,N}-\frac{1}{2}][1]}
{[v_2-v'-1]^*[P_{N-1,N}-1]^*
[v_1-v'+\frac{1}{2}][P_{N-1,N}+h_{N-1,N}]}
\right.\nonumber\\
&-&\left.
\oint_{{C}_-} \frac{dz'}{2\pi i z'} 
H_{N-1}^-\left(v'-\frac{1}{4}\right)
\frac{[v_2-v'-P_{N-1,N}+1]^*[1]^*
[v_1-v'+P_{N-1,N}+
h_{N-1,N}-\frac{1}{2}][1]}
{[v_2-v']^*[P_{N-1,N}-1]^*
[v_1-v'+\frac{1}{2}][P_{N-1,N}+h_{N-1,N}]}
\right).\nonumber\\
&&\label{proof:Type-I,II}
\end{eqnarray}
The contours are specified by
\begin{eqnarray}
&&{C}_+:|qz_1|,|q^{-2}z_2|<|z'|<|p^{-1}qz_1|,
|p^{*-1}q^{-2}z_2|,
|p^{-1}q^{-1}z_1|, |p^{*-1}z_2|,\\
&&{C}_-:|qz_1|,|z_2|<|z'|<|p^{-1}qz_1|,
|p^{*-1}z_2|, |q^{-1}z_1|, |q^2z_2|.
\end{eqnarray}
Here 
the conditions $|z'|<
|p^{-1}q^{-1}z_1|,|p^{*-1}z_2|$ for $C_+$
 and $|z'|<|q^{-1}z_1|,|q^2z_2|$ for $C_-$ are added
because of the convergence of the operator 
product $\Phi_N(v_1)\Psi_N^*(v_2)H_{N-1}^+(v'+1/4)$
and $\Phi_N(v_1)\Psi_N^*(v_2)H_{N-1}^-(v'-1/4)$, respectively.
Changing the integration variable $z' \to pz'$ in the second term
and using the periodicity of $[v], [v]^*$ and the relation 
$H_{N-1}^-(v'+r-\frac{1}{4})=H_{N-1}^+(v'+\frac{1}{4})$, 
we see that the integrand in the second term coincides with 
the one in the first
term but the contour in the second term is changed to
\be
&&\tilde{C}_-:|p^{-1}qz_1|,|p^{*-1}q^{-2}z_2|<|z'|<
|p^{-2}qz_1|,
|p^{*-2}q^{-2}z_2|, 
|p^{-1}q^{-1}z_1|,
|p^{*-1}z_2|.
\en
Here $\tilde{C}_-$ encircles the same poles as $C_+$. In addition,
$\tilde{C}_-$ would encircle 
two extra poles at $z'=p^{-1}qz_1,\ p^{*-1}q^{-2}z_2$, if 
the operator product $\Phi_N(v_1)\Psi_N^*(v_2)
H_{N-1}^+(v'+\frac{1}{4})$ had no zeros which cancel these extra 
poles. In fact, the operator product does have zeros at the 
required points. Therefore the RHS of \eqref{proof:Type-I,II} vanishes.
The proof for the general $1\leq j_1,  j_2 \leq N-1$ case is similar.
\begin{flushright}
{Q.E.D.}
\end{flushright}

Now let us investigate the intertwining relation 
for level $c=1$.
For this purpose, we remind the reader of the fact that 
in the trigonometric case, i.e. $U_q(\slnh)$,  
the $L$-operator can be constructed as a composition 
of type I and II vertex operators \cite{Miki,FIJKMY}.
The following theorem
is an elliptic analogue of such a construction.

\begin{thm}~~~For $c=1$, the $L$-operator $\widehat{L}^+(v)$
is given by a product of the type-I and type-II
vertex operators.
\begin{eqnarray}
\widehat{L}^+_{j k}(v)
=\frac{1}{g_N}\ 
\Psi_{k}^*(v+r)
\Phi_{j}\left(v+r+\frac{1}{2}\right)
~~\qquad (1\leq j,k \leq N).
\label{L-vertex}
\end{eqnarray}
Here we set
\begin{eqnarray}
g_N=\frac{(q^{2};q^{2N})_\infty}
{(q^{2N};q^{2N})_\infty}.
\end{eqnarray}
\end{thm}
{\it Proof.}~~~
For the special component $j=k=N$ of
(\ref{L-vertex}), we have 
\bea
&&K^+_N(v)=g_N^{-1}
\Psi_{N}^*(v+r)
\Phi_{N}\left(v+r+\frac{1}{2}\right). \lb{kpsiphi}
\ena
This is a direct consequence of 
\eqref{free:Phi} and \eqref{free:Psi}.
Let us consider
the $j<k=N$ component of (\ref{L-vertex}).
After a few calculations using \eqref{VII2} and \eqref{kpsiphi},
we can reduce this  to relation (\ref{VI4}).
Similarly, the $N=j>k$ component of (\ref{L-vertex})
is reduced to
 relation (\ref{VII4}).

Next, let us study
the simplest non-trivial component $j=k=N-1$.
From \eqref{VI4}, \eqref{VII2}-\eqref{VII6}, \eqref{k5}-\eqref{k6} 
and \eqref{EiFj}, 
we have the following equation after integrating
the delta functions.
\begin{eqnarray}
&&g_N^{-1}\Psi_{N-1}^*(v+r)\Phi_{N-1}(v+r+1/2)-
F_{N-1,N}^+(v)K_N^+(v)E_{N,N-1}^+(v)\nonumber\\
&=&\frac{a_{N-1,N}a_{N-1,N}^*}{q-q^{-1}}\nn\\
&\times&
\left(
\oint_{{C}_+} \frac{dz'}
{2\pi i z'} 
H_{N-1}^+\left(v'+\frac{1}{4}
\right)K_N^+(v)
\frac{[v-v'-P_{N-1,N}+1]^*[1]^*
[v-v'+P_{N-1,N}+h_{N-1,N}][1]}
{[v-v'+1]^*[P_{N-1,N}-1]^*
[v-v'][P_{N-1,N}+h_{N-1,N}]}
\right.\nonumber\\
&-&\left.
\oint_{{C}_-} \frac{dz'}{2\pi i z'} 
H_{N-1}^-\left(v'-\frac{1}{4}\right)
K_N^+(v)
\frac{[v-v'-P_{N-1,N}+2]^*[1]^*
[v-v'+P_{N-1,N}+h_{N-1,N}][1]}
{[v-v'+2]^*[P_{N-1,N}-1]^*
[v-v'][P_{N-1,N}+h_{N-1,N}]}
\right).\nonumber\\
&&\label{proof:L=TypeI,II}
\end{eqnarray}
Here the contours are
\begin{eqnarray}
&&{C}_+:|pz|,|p^*z|<|z'|<|z|,\\
&&{C}_-:|pz|<|z'|<|z|,|q^2z|.
\end{eqnarray}
Changing the integration variable $z' \to pz'$
in the second term and using the periodicity of
$[v],[v]^*$ and the relation $H_{N-1}^-(v'+r-1/4)=H_{N-1}^+(v'+1/4)$,
we see that the integrand in the second term coincides with 
the first term, but the contour in the second term is changed
to
\begin{eqnarray}
\tilde{C}_-:|z|<|z'|<|p^{-1}z|,|p^{*-1}z|.
\end{eqnarray}
The contour $\tilde{C}_-$ encircles
the same poles as $C_+$ together with one additional pole
at $z'=z$. Hence the RHS of \eqref{proof:L=TypeI,II} 
becomes the residue at $z'=z$. We thus obtain
\begin{eqnarray}
g_N^{-1}\Psi_{N-1}^*(v+r)\Phi_{N-1}(v+r+1/2)=
F_{N-1,N}^+(v)K_N^+(v)E_{N,N-1}^+(v)+K_{N-1}^+(v).
\end{eqnarray}
The RHS coincides with the $(N-1,N-1)$ component of $\widehat{L}^+(v)$.
The proof of the general case $1\leq j,k \leq N-1$ 
is similar.
\begin{flushright}
{Q.E.D.}
\end{flushright}

\begin{cor}~~~
For $c=1$, the $L$-operator $\widehat{L}^+(v)$ satisfies
the $RLL=LLR^*$ relation \eqref{thm:RLL}.
\end{cor}
{\it Proof.}~~~
Let us substitute the expression \eqref{L-vertex} into the LHS of
\eqref{thm:RLL}. Then using 
the commutation relations of the vertex operators
(\ref{Com:Type-I})-(\ref{Com:Type-I,II}) and the formula
\begin{eqnarray}
\frac{\rho^+(v)}{
\rho^{+*}(v)
}=\frac{\tilde{\rho}^+(v)}{\tilde{\rho}^{+*}(v)}
\frac{{\chi}(\frac{1}{2}-v)}{
{\chi}(\frac{1}{2}+v)},
\end{eqnarray}
one gets the desired result.
\begin{flushright}
{Q.E.D.}
\end{flushright}

In the same way, we have 
\begin{cor}
~~~For 
$c=1$, the type-I and the type II vertex operators $\widehat{\Phi}_V(v)$,
 $\widehat{\Psi}^*_V(v)$
 satisfy the full intertwining relations \eqref{Inter:Uqp1} and 
 \eqref{Inter:Uqp2} 
 with $V=W\cong \C^N$.
\end{cor}



~\\
{\bf Acknowledgements}~~
The authors would like to thank Robert Weston and colleagues in 
the Department of Mathematics,
Heriot-Watt University, where a part of this work was done, for 
their kind hospitality.
H.K is also grateful to JSPS and the Royal Society for the exchange fellowship.
This work is also supported by 
Grant-in-Aid for Scientific Research
({\bf C}) (11640030, 14540028)
and Grant-in-Aid for Young Scientist
({\bf B}) (14740107) from
the Ministry of Education, Science, Sports and Culture.

\begin{appendix}

\section{Operator Product Expansions}\lb{OPE}
Here we list formulae of operator product expansions (OPE) used in 
Section 3.3 and 6.2. For operators $A(z),\ B(z)$, we write
\be
&&A(z_1)B(z_2)=\lal A(z_1)B(z_2)\ral :A(z_1)B(z_2):.
\en

\noindent
(I) In Section 3.3, we used the OPE's of the currents $\psi_j^\pm(z,p)$ (\ref{def:psi+})-(\ref{def:psi-})
 and $k_j(z,p)$ (\ref{def:k}) 
for generic $c$:
\begin{eqnarray}
&&\lal k_j(z_1,p)k_j(z_2,p)\ral=
\frac{\{pq^2z_2/z_1\} 
\{pq^{2N-2}z_2/z_1\}
\{p^*q^{2N}z_2/z_1\}^* 
\{p^*z_2/z_1\}^*
}
{\{pq^{2N}z_2/z_1\} 
\{pz_2/z_1\}
\{p^*q^2z_2/z_1\}^* 
\{p^*q^{2N-2}z_2/z_1\}^*
},\\
&&\lal k_{j_1}(z_1,p)k_{j_2}(z_2,p)\ral=
\frac{\{pq^{2N+2}z_2/z_1\} \{pq^{2N-2}z_2/z_1\}
\{p^*q^{2N}z_2/z_1\}^{*2}
}{
\{pq^{2N}z_2/z_1\}^2
\{p^*q^{2N+2}z_2/z_1\}^* \{p^*q^{2N-2}z_2/z_1\}^*
}
\qquad (j_1<j_2),\nonumber\\
\\&&
\lal k_{j_2}(z_1,p)k_{j_1}(z_2,p)\ral=
\frac{\{pq^{2}z_2/z_1\} \{pq^{-2}z_2/z_1\}
\{p^* z_2/z_1\}^{*2}
}{
\{p z_2/z_1\}^2
\{p^*q^{2}z_2/z_1\}^* 
\{p^*q^{-2}z_2/z_1\}^*
}
\qquad (j_1<j_2).
\end{eqnarray}
\begin{eqnarray}
&&\lal\psi_j^+(z_1,p)\psi_j^+(z_2,p)\ral=\frac{
(pq^2z_2/z_1;p)_\infty
(p^*q^{-2}z_2/z_1;p^*)_\infty
}{
(pq^{-2}z_2/z_1;p)_\infty
(p^*q^2 z_2/z_1;p^*)_\infty
},\\&&
\lal\psi_j^+(z_1,p)\psi_{j+1}^+(z_2,p)\ral
=\frac{(pq^{-1}z_2/z_1;p)_\infty
(p^*q z_2/z_1;p^*)_\infty
}
{
(pqz_2/z_1;p)_\infty
(p^*q^{-1}z_2/z_1;p^*)_\infty
},\\&&
\lal\psi_{j+1}^+(z_1,p)\psi_j^+(z_2,p)\ral=
\frac{(pq^{-1}z_2/z_1;p)_\infty
(p^*q z_2/z_1;p^*)_\infty
}
{
(pqz_2/z_1;p)_\infty
(p^*q^{-1}z_2/z_1;p^*)_\infty
}.
\end{eqnarray}

\noindent
(II) In Section 6.2, we used the OPE's of 
the currents $\Phi_N(v),\ \Psi^*_N(v),\ E_j(v),\ F_j(v),\ K^+_j(v)$
 for $c=1$. We here list their boson part only. Namely, let us define the 
 boson part of them by
\begin{eqnarray}
&&\phi_N(v)=:\exp\left(-\sum_{m \neq 0}\frac{1}{m}B_m^N
z^{-m}\right):,\\
&&\psi_N^*(v)=:\exp\left(
\sum_{m \neq 0}\frac{[rm]_q}{m[r^*m]_q}B_m^N z^{-m}
\right):,\\
&&e_j(v)=:\exp\left(-\sum_{m \neq 0}\frac{[rm]_q}{m[r^*m]_q}
(-B_m^j+B_m^{j+1})(q^{N-j}z)^{-m}\right):,\\
&&f_j(v)=:\exp\left(
\sum_{m\neq 0}\frac{1}{m}(-B_m^j+B_m^{j+1})(q^{N-j}z)^{-m}
\right):.
\end{eqnarray}
Then the OPE's of them are given by
\begin{eqnarray}
&&\langle K_j^+(v_1)f_j(v_2)\rangle=
\frac{(q^{N-j+1}z_2/z_1;p)_\infty}{
(q^{N-j-1}z_2/z_1;p)_\infty
},\\
&&\langle K_{j+1}^+(v_1)f_j(v_2)\rangle=
\frac{(q^{N-j-3}z_2/z_1;p)_\infty}{
(q^{N-j-1}z_2/z_1;p)_\infty},\\
&&\langle f_j(v_1)K_j^+(v_2)\rangle=
\frac{(q^{-N+j-1}z_2/z_1;p)_\infty}{
(q^{-N+j-3}z_2/z_1;p)_\infty
},\\
&&\langle f_j(v_1)K_{j+1}^+(v_2)\rangle=
\frac{(q^{-N+j-1}z_2/z_1;p)_\infty}{
(q^{-N+j+1}z_2/z_1;p)_\infty},
\end{eqnarray}
\begin{eqnarray}
&&\langle K_j^+(v_1)e_j(v_2)\rangle=
\frac{(q^{N-j-2}z_2/z_1;p^*)_\infty}{
(q^{N-j}z_2/z_1;p^*)_\infty
},\\
&&\langle K_{j+1}^+(v_1)e_j(v_2)\rangle=
\frac{(q^{N-j-2}z_2/z_1;p^*)_\infty}{
(q^{N-j-4}z_2/z_1;p^*)_\infty},\\
&&\langle e_j(v_1)K_j^+(v_2)\rangle=
\frac{(q^{-N+j-4}z_2/z_1;p^*)_\infty}{
(q^{-N+j-2}z_2/z_1;p^*)_\infty
},\\
&&\langle e_j(v_1)K_{j+1}^+(v_2)\rangle=
\frac{(q^{-N+j}z_2/z_1;p^*)_\infty}{
(q^{-N+j-2}z_2/z_1;p^*)_\infty},
\end{eqnarray}
\begin{eqnarray}
&&\lal\phi_N(v_1)K_N^+(v_2)\ral=
\frac{\{pq^3z_2/z_1\} \{pq^{2N-1}z_2/z_1\}}
{\{pqz_2/z_1\} \{pq^{2N+1}z_2/z_1\}},\\&&
\lal K_N^+(v_2)\phi_N(v_1)\ral=
\frac{\{qz_1/z_2\}\{q^{2N-3}z_1/z_2\}}
{\{q^{-1}z_1/z_2\} \{q^{2N-1}z_1/z_2\}},\\&&
\lal\phi_N(v_1)K_j^+(v_2)\ral=
\frac{\{pq^3z_2/z_1\} 
\{pq^{-1}z_2/z_1\}}{
\{pqz_2/z_1\}^2
},\\&&
\lal K_j^+(v_2)\phi_N(v_1)\ral=
\frac{\{q^{2N+1}z_1/z_2\}
\{q^{2N-3}z_1/z_2\}
}{\{q^{2N-1}z_1/z_2\}^2},\\
&&\nn\\&&
\lal\psi_N^*(v_1)K_N^+(v_2)\ral=
\frac{
\{p^*q^2 z_2/z_1\}^* 
\{p^*q^{2N+2}z_2/z_1\}^*}
{\{p^*q^4z_2/z_1\}^* \{p^*q^{2N}z_2/z_1\}^*}
,\\&&
\lal K_N^+(v_2)\psi_N^*(v_1)\ral=
\frac{\{q^{-2}z_1/z_2\}^* \{q^{2N-2}z_1/z_2\}^*}
{\{z_1/z_2\}^* \{q^{2N-4}z_1/z_2\}^*},\\&&
\lal\psi_N^*(v_1)K_j^+(v_2)\ral=
\frac{\{p^*q^2z_2/z_1\}^{*2}}{
\{p^*q^4z_2/z_1\}^* \{p^*z_2/z_1\}^*
},\\&&
\lal K_j^+(v_2)\psi_N^*(v_1)\ral=
\frac{\{q^{2N-2}z_1/z_2\}^{*2}}{
\{q^{2N}z_1/z_2\}^* \{q^{2N-4}z_1/z_2\}^*
},\\
&&\nn\\
&&\lal f_{N-1}(v_1)\phi_N(v_2)\ral=
\frac{(pq^{-1}z_2/z_1;p)_\infty}{(qz_2/z_1;p)_\infty},\\&&
\lal \phi_N(v_1)f_{N-1}(v_2)\ral=
\frac{(pq^{-1}z_2/z_1;p)_\infty}{(qz_2/z_1;p)_\infty},\\&&
\lal e_{N-1}(v_1)\psi_N^*(v_2)\ral=
\frac{(p^{*}qz_2/z_1;p^*)_\infty}
{(q^{-1}z_2/z_1;p^*)_\infty},\\&&
\lal\psi_N^*(v_1)e_{N-1}(v_2)\ral=
\frac{(p^*qz_2/z_1;p^*)_\infty}
{(q^{-1} z_2/z_1;p^*)_\infty},\\
&&\nn\\
&&\lal\phi_N(v_1)\phi_N(v_2)\ral=\frac{\{pq^{2N-2}z_2/z_1\}
\{q^2 z_2/z_1\}}{
\{q^{2N}z_2/z_1\}
\{pz_2/z_1\}},\\&&
\lal\psi_N^*(v_2)\psi_N^*(v_1)\ral=
\frac{\{p^*q^{2N}z_2/z_1\}^* 
\{z_2/z_1\}^*}{
\{p^*q^{2}z_2/z_1\}^* 
\{q^{2N-2}z_2/z_1\}^*
},\\&&
\lal\phi_N(v_1)\psi_N^*(v_2)\ral=
\frac{(q^{2N-1}z_2/z_1;q^{2N})_\infty }
{(qz_2/z_1;q^{2N})_\infty},\\&&
\lal\psi_N^*(v_1)\phi_N(v_2)\ral=
\frac{(q^{2N-1}z_2/z_1;q^{2N})_\infty }{(qz_2/z_1;q^{2N})_\infty},\\
&&\nn\\
&&\lal e_j(v_1)e_j(v_2)\ral=
\frac{(z_2/z_1;p^*)_\infty 
(q^{-2}z_2/z_1;p^*)_\infty}
{(p^*q^2z_2/z_1;p^*)_\infty 
(p^*z_2/z_1;p^*)_\infty}
,\\&&
\lal e_j(v_1)e_{j+1}(v_2)\ral=
\frac{(p^*qz_2/z_1;p^*)_\infty}
{(q^{-1}z_2/z_1;p^*)_\infty},\\&&
\lal e_{j+1}(v_1)e_j(v_2)\ral=
\frac{(p^*qz_2/z_1;p^*)_\infty}
{(q^{-1}z_2/z_1;p^*)_\infty}
,\\&&
\lal f_j(v_1)f_j(v_2)\ral=
\frac{(z_2/z_1;p)_\infty (q^2z_2/z_1;p)_\infty}{
(pz_2/z_1;p)_\infty (pq^{-2}z_2/z_1;p)_\infty}
,\\&&
\lal f_j(v_1)f_{j+1}(v_2)\ral=
\frac{(pq^{-1}z_2/z_1;p)_\infty}{
(qz_2/z_1;p)_\infty}
,\\&&
\lal f_{j+1}(v_1)f_j(v_2)\ral=
\frac{(pq^{-1}z_2/z_1;p)_\infty}{
(qz_2/z_1;p)_\infty}.
\end{eqnarray}


\section{Proof of the Relations \eqref{eljelj}-\eqref{fjlfjl} and 
\eqref{elkelj}-\eqref{fklfjl}}\label{Proof}
Let us consider the relations
\bea
&&K^+_{l}(v_2)^{-1}E^+_{l,j}(v_1)K_l^+(v_2)E^+_{l,j}(v_2)=K^+_{l}(v_1)^{-1}
E^+_{l,j}(v_2)K_l^+(v_1)E^+_{l,j}(v_1),\lb{eljelj2}\\
&&F^+_{j,l}(v_1)K^+_{l}(v_1)F^+_{j,l}(v_2)K_l^+(v_1)^{-1}=
F^+_{j,l}(v_2)K^+_{l}(v_2)F^+_{j,l}(v_1)K_l^+(v_2)^{-1},\lb{fjlfjl2}
\ena
for $1\leq j<l\leq N$. Then the relations \eqref{eljelj} and \eqref{fjlfjl} 
follow from these relations and
\eqref{klelj}, \eqref{klfjl}.
In this Appendix, we give proofs of the relations \eqref{fjlfjl2} 
and \eqref{fklfjl}. 
The proof of the other cases \eqref{eljelj2} and \eqref{elkelj} are similar.

Let us set 
\bea
&&f(v,w)=\frac{[v+\frac{1}{2}-w]}{[v-\frac{1}{2}]},\\
&&h(v)=\frac{[v-1]}{[v+1]}.
\ena
Recall that the half current $F^+_{j,l}(v)$ is given by
\begin{eqnarray}
F_{j,l}^+(v)&=&a_{j,l}\oint_{C(j,l)}
\prod_{m=j}^{l-1}\frac{dz_k}{2\pi i z_k}
F_{l-1}(v_{l-1})F_{l-2}(v_{l-2})
\cdots F_{j}(v_{j})\nonumber\\
&&\times
\prod_{m=j}^{l-1}f(v_{m}-v_{m+1},\pi_{m+1,j})\frac{[1]}{[\pi_{m+1,j}-\delta_{m,l-1}]}
,\lb{fjl2}
\end{eqnarray}
where we set $v_l=v+\frac{l-N-1}{2}$. Recall also 
$z_m=q^{2v_m}$ and $\pi_{l,j}=P_{j,l}+h_{j,l}$.
We call $F_{l-1}(v_{l-1})$ $F_{l-2}(v_{l-2})
\cdots F_{j}(v_{j})$ the operator part, and $\prod_{m=j}^{l-1}f(v_{m}-v_{m+1},
\pi_{m+1,j})\frac{[1]}{[\pi_{m+1,j}-\delta_{m,l-1}]}$ the coefficient part.
We keep coefficient parts in the right of operator parts.
In the coefficient part, We represent $\prod_{m=j}^{\l-1}
f(v_{m}-v_{m+1},\pi_{m+1,j})$ by the diagram
\be
\begin{diagram}
\node{v_{j}}\arrow{e,t}{\pi_{j+1,j}}
\node{v_{j+1}}\arrow{e,t}{\pi_{j+2,j}}
\node{\phantom{f}\cdots\phantom{f}}\arrow{e,t}{\pi_{l-1,j}}
\node{\phantom{f}v_{l-1}\phantom{f}}\arrow{e,t}{\pi_{l,j}}
\node{\phantom{f}v_{l}.\phantom{f}}
\end{diagram}
\en
According to the relation \eqref{FiFj} with $i=j$, we have the equality

\be
\oint
{dz_j\over2\pi iz_j}{dz'_j\over2\pi iz'_j}
F_j(v_j)F_j(v'_j)A(v_j,v'_j)
=\oint
{dz_j\over2\pi iz_j}{dz'_j\over2\pi iz'_j}
F_j(v_j)F_j(v'_j)h(v'_j-v_j)A(v'_j,v_j),
\en
when the integration contours for $z_j$ and $z'_j$ are the same.
We define `weak equality' in the following sense\cite{AJMP}. 
The two coefficient functions $A(v_j,v'_j)$ and $B(v_j,v'_j)$ coupled to
$F_j(v_j)F_j(v'_j)$ in integrals are equal in weak sense if
\be
A(v_j,v'_j)+h(v'_j-v_j)A(v'_j,v_j)
=B(v_j,v'_j)+h(v'_j-v_j)B(v'_j,v_j).
\en
We write the weak equality as
\be
A(v_j,v'_j)\sim B(v_j,v'_j).
\en

To prove the equality \eqref{fjlfjl2} and \eqref{fklfjl},
it is enough to show the equalities of coefficient parts in weak sense.

Let us recall the following two lemmas\cite{AJMP}.
\begin{lem}\lb{lem1}
The coefficient function
\bea
\begin{diagram}
\node{v_{j}}\arrow{e,t}{\pi_{j+1,j}-1}\arrow{se,t}{0}
\node{v_{j+1}}\arrow{e,t}{\pi_{j+2,j}-1}\arrow{se,t}{0}
\node{\phantom{f}\cdots\phantom{f}}\arrow{e,t}{\pi_{l-2,j}-1}\arrow{se,t}{0}
\node{\phantom{f}v_{l-2}\phantom{f}}\arrow{e,t}{\pi_{l-1,j}-1}\arrow{se,t}{0}
\node{\phantom{f}v_{l-1}\phantom{f}}\\
\node{v'_{j}}\arrow{e,t}{\pi_{j+1,j}}
\node{v'_{j+1}}\arrow{e,t}{\pi_{j+2,j}}
\node{\phantom{f}\cdots\phantom{f}}\arrow{e,t}{\pi_{l-2,\mu}}
\node{\phantom{f}v'_{l-2}\phantom{f}}\arrow{e,t}{\pi_{l-1,\mu}}
\node{\phantom{f}v'_{l-1}\phantom{f}}
\end{diagram}
\ena
is invariant in weak sense when $v_{l-1}$ and $v'_{l-1}$ are exchanged. 
\end{lem}

\begin{lem}\lb{lem2}
\bea
&&\begin{diagram}
\node{{v'_{k-1}}}\arrow{e,t}{\pi_{k,j}+1}\arrow{se,r}{0}
\node{v_{k}}\arrow{e,t}{\pi_{k+1,j}}\arrow{se,r}{0}
\node{v_{k+1}}\arrow{e,t}{\pi_{k+2,j}}\arrow{se,r}{0}
\node{\cdots}\arrow{e,t}{\pi_{l-2,j}}\arrow{se,r}{0}
\node{v_{l-2}}\arrow{e,t}{\pi_{l-1,j}}\arrow{se,r}{0}
\node{v_{l-1}}\\
\node[2]{\phantom{f}v'_{k}\phantom{f}}\arrow{e,b}{\pi_{k+1,k}}
\node{v'_{k+1}}\arrow{e,b}{\pi_{k+2,k}}
\node{\cdots}\arrow{e,b}{\pi_{l-2,k}}
\node{v'_{l-2}}\arrow{e,b}{\pi_{l-1,k}}
\node{v'_{l-1}}
\end{diagram}\nonumber\\
&&\sim\frac{1}{\beta(v_{l-1}-v'_{l-1},\pi_{k,j})}\nn\\
&&\qquad\times\begin{diagram}
\node{{v'_{k-1}}}\arrow{e,t}{\pi_{k,j}}
\node{v'_{k}}\arrow{e,t}{\pi_{k+1,j}}
\node{v'_{k+1}}\arrow{e,t}{\pi_{k+2,j}}
\node{\phantom{f}\cdots}\arrow{e,t}{\pi_{l-2,j}}
\node{v'_{l-2}}\arrow{e,t}{\pi_{l-1,j}}
\node{v_{l-1}}\\
\node[2]{\phantom{f}v_{k}\phantom{f}}\arrow{e,b}{\pi_{k+1,k}}\arrow{ne,r}{0}
\node{v_{k+1}}\arrow{e,b}{\pi_{k+2,k}}\arrow{ne,r}{0}
\node{\cdots}\arrow{e,b}{\pi_{l-2,k}}\arrow{ne,r}{0}
\node{v_{l-2}}\arrow{e,b}{\pi_{l-1,k}}\arrow{ne,r}{0}
\node{v'_{l-1}}
\end{diagram}\nonumber\\
&&-\frac{\gamma(v_{l-1}-v'_{l-1},\pi_{k,j})}{\beta(v_{l-1}-v'_{l-1},\pi_{k,j})}\nonumber\\
&&\qquad\times\begin{diagram}
\node{{v'_{k-1}}}\arrow{e,t}{\pi_{k,j}}
\node{v'_{k}}\arrow{e,t}{\pi_{k+1,j}}
\node{v'_{k+1}}\arrow{e,t}{\pi_{k+2,j}}
\node{\phantom{f}\cdots}\arrow{e,t}{\pi_{l-2,j}}
\node{v'_{l-2}}\arrow{e,t}{\pi_{l-1,j}}
\node{v'_{l-1}}\\
\node[2]{\phantom{f}v_{k}\phantom{f}}\arrow{e,b}{\pi_{k+1,k}}\arrow{ne,r}{0}
\node{v_{k+1}}\arrow{e,b}{\pi_{k+2,k}}\arrow{ne,r}{0}
\node{\cdots}\arrow{e,b}{\pi_{l-2,k}}\arrow{ne,r}{0}
\node{v_{l-2}}\arrow{e,b}{\pi_{l-1,k}}\arrow{ne,r}{0}
\node{v_{l-1}.}
\end{diagram}\nonumber\\
\ena
where
\be
\beta(v,w)={[v][w-1]\over[v+1][w]},\quad \gamma(v,w)={[v+w][1]\over[v+1][w]}.
\en
\end{lem}

Now let us  show the relation \eqref{fjlfjl2}. 
By using \eqref{k7}-\eqref{k9}, \eqref{FiFj}, 
\eqref{kjpjl} and \eqref{fjpjlphjl} in the LHS of \eqref{fjlfjl2}, the 
operator part can be arranged 
to $F_{l-1}(v_{l-1})F_{l-1}(v'_{l-1})F_{l-2}(v_{l-2})F_{l-2}(v'_{l-2})
\cdots F_{j}(v_{j})F_{j}(v'_{j})$. Then the coefficient part is given by 
the product of the  factors 
$\prod_{m=j}^{l-1}\frac{[1]^2}{[\pi_{m+1,j}-1-2\delta_{k,l-1}]
[\pi_{m+1,j}-2\delta_{k,l-1}]}$ and the one
represented by the diagram
\bea
\begin{diagram}
\node{v_{j}}\arrow{e,t}{\pi_{j+1,j}-1}\arrow{se,t}{0}
\node{v_{j+1}}\arrow{e,t}{\pi_{j+2,j}-1}\arrow{se,t}{0}
\node{\phantom{f}\cdots\phantom{f}}\arrow{e,t}{\pi_{l-2,j}-1}\arrow{se,t}{0}
\node{\phantom{f}v_{l-2}\phantom{f}}\arrow{e,t}{\pi_{l-1,j}-1}\arrow{se,t}{0}
\node{\phantom{f}v_{l-1}\phantom{f}}\arrow{e,t}{\pi_{l,j}-2}
\node{\phantom{f}v_{l}\phantom{f}}\\
\node{v'_{j}}\arrow{e,t}{\pi_{j+1,j}}
\node{v'_{j+1}}\arrow{e,t}{\pi_{j+2,j}}
\node{\phantom{f}\cdots\phantom{f}}\arrow{e,t}{\pi_{l-2,\mu}}
\node{\phantom{f}v'_{l-2}\phantom{f}}\arrow{e,t}{\pi_{l-1,\mu}}
\node{\phantom{f}v'_{l-1}\phantom{f}}\arrow{e,t}{\pi_{l,\mu}-1}\arrow{ne,t}{2}
\node{\phantom{f}v'_{l}.\phantom{f}}
\end{diagram} \lb{b8}
\ena
The relation \eqref{fjlfjl2} denotes that \eqref{b8} is invariant, 
at least in weak sense, when $v_l$ and $v'_l$ are exchanged. 
Applying the Lemma \ref{lem1} to the corresponding part of \eqref{b8}, 
it is enough to show the weak equality for the rest part
\bea
&&f(v_{l-1}-v_l,\pi_{l,j}-2)f(v'_{l-1}-v'_l,\pi_{l,j}-1)f(v'_{l-1}-v_l,2)\nn\\
&&\qquad\qquad\sim 
f(v_{l-1}-v'_l,\pi_{l,j}-2)f(v'_{l-1}-v_l,\pi_{l,j}-1)f(v'_{l-1}-v'_l,2).
\lb{fjlfjl3}
\ena
Let us set $v=v_{l-1}, v'=v_{l-1}, w=\pi_{l,j}$ and denote the LHS and the RHS 
by $A(v,v')$ and $B(v,v')$, respectively. Then from the theta function identity
such as \eqref{thetaid}, we have 
\be
A(v,v')-B(v,v')=\frac{[v-v'+1][v+v'-v_l-v'_l-w][v_l-v'_l][w-2]}
{[v-v_l-\frac{1}{2}][v'-v'_l-\frac{1}{2}][v-v'_l-\frac{1}{2}][v'-v_l-\frac{1}{2}]}.
\en
Then it is easy to show
\be
h(v'-v)(A(v',v)-B(v',v))=-(A(v,v')-B(v,v')). 
\en
Therefore we get the weak equality \eqref{fjlfjl3}.

Next we prove \eqref{fklfjl} ($j<k<l$). 
The equality follows from the weak equality (A)+(B)+(C)$\sim$0,
where
\bea
&&(A)=\nn\\
&&\begin{diagram}
\node{{v'_{k-1}}}\arrow{e,t}{\pi_{k,j}}
\node{v'_{k}}\arrow{e,t}{\pi_{k+1,j}}
\node{v'_{k+1}}\arrow{e,t}{\pi_{k+2,j}}
\node{\phantom{f}\cdots}\arrow{e,t}{\pi_{l-2,j}}
\node{v'_{l-2}}\arrow{e,t}{\pi_{l-1,j}}
\node{v'_{l-1}}\arrow{e,t}{\pi_{l,j}-1}\arrow{se,r}{2}
\node{v'_l}\\
\node[2]{\phantom{f}v_{k}\phantom{f}}\arrow{e,b}{\pi_{k+1,k}}\arrow{ne,r}{0}
\node{v_{k+1}}\arrow{e,b}{\pi_{k+2,k}}\arrow{ne,r}{0}
\node{\cdots}\arrow{e,b}{\pi_{l-2,k}}\arrow{ne,r}{0}
\node{v_{l-2}}\arrow{e,b}{\pi_{l-1,k}}\arrow{ne,r}{0}
\node{v_{l-1}}\arrow{e,b}{\pi_{l,k}-1}
\node{v_{l},}
\end{diagram}\nonumber\\
&&(B)=-b(v'_l-v_l,\pi_{k,j})\frac{[\pi_{k,j}]}{[\pi_{k,j}+1]}
\times\nonumber\\
&&\begin{diagram}
\node{{v'_{k-1}}}\arrow{e,t}{\pi_{k,j}+1}\arrow{se,r}{0}
\node{v_{k}}\arrow{e,t}{\pi_{k+1,j}}\arrow{se,r}{0}
\node{v_{k+1}}\arrow{e,t}{\pi_{k+2,j}}\arrow{se,r}{0}
\node{\cdots}\arrow{e,t}{\pi_{l-2,j}}\arrow{se,r}{0}
\node{v_{l-2}}\arrow{e,t}{\pi_{l-1,j}}\arrow{se,r}{0}
\node{v_{l-1}}\arrow{e,t}{\pi_{l,j}-1}
\node{v'_l}\\
\node[2]{\phantom{f}v'_{k}\phantom{f}}\arrow{e,b}{\pi_{k+1,k}}
\node{v'_{k+1}}\arrow{e,b}{\pi_{k+2,k}}
\node{\cdots}\arrow{e,b}{\pi_{l-2,k}}
\node{v'_{l-2}}\arrow{e,b}{\pi_{l-1,k}}
\node{v'_{l-1}}\arrow{e,b}{\pi_{l,k}-1}\arrow{ne,r}{2}
\node{v_{l},}
\end{diagram}\nonumber\\
&&(C)=-c(v'_l-v_l,\pi_{k,j})\times\nonumber\\
&&\begin{diagram}
\node{{v'_{k-1}}}\arrow{e,t}{\pi_{k,j}}
\node{v'_{k}}\arrow{e,t}{\pi_{k+1,j}}
\node{v'_{k+1}}\arrow{e,t}{\pi_{k+2,j}}
\node{\phantom{f}\cdots}\arrow{e,t}{\pi_{l-2,j}}
\node{v'_{l-2}}\arrow{e,t}{\pi_{l-1,j}}
\node{v'_{l-1}}\arrow{e,t}{\pi_{l,j}-1}\arrow{se,r}{2}
\node{v_l}\\
\node[2]{\phantom{f}v_{k}\phantom{f}}\arrow{e,b}{\pi_{k+1,k}}\arrow{ne,r}{0}
\node{v_{k+1}}\arrow{e,b}{\pi_{k+2,k}}\arrow{ne,r}{0}
\node{\cdots}\arrow{e,b}{\pi_{l-2,k}}\arrow{ne,r}{0}
\node{v_{l-2}}\arrow{e,b}{\pi_{l-1,k}}\arrow{ne,r}{0}
\node{v_{l-1}}\arrow{e,b}{\pi_{l,k}-1}
\node{v'_{l}.}
\end{diagram}\nonumber\\
\ena

Using the weak equality in Lemma \ref{lem2}, we modify $(B)$ to $(A')+(C')$ where
\bea
&&(A')=-\frac{b(v'_l-v_l,\pi_{k,j})}{\beta(v_{l-1}-v'_{l-1},\pi_{k,j})}\frac{[\pi_{k,j}]}{[\pi_{k,j}+1]}
\times\nonumber\\
&&\begin{diagram}
\node{{v'_{k-1}}}\arrow{e,t}{\pi_{k,j}}
\node{v'_{k}}\arrow{e,t}{\pi_{k+1,j}}
\node{v'_{k+1}}\arrow{e,t}{\pi_{k+2,j}}
\node{\phantom{f}\cdots}\arrow{e,t}{\pi_{l-2,j}}
\node{v'_{l-2}}\arrow{e,t}{\pi_{l-1,j}}
\node{v_{l-1}}\arrow{e,t}{\pi_{l,j}-1}
\node{v'_l}\\
\node[2]{\phantom{f}v_{k}\phantom{f}}\arrow{e,b}{\pi_{k+1,k}}\arrow{ne,r}{0}
\node{v_{k+1}}\arrow{e,b}{\pi_{k+2,k}}\arrow{ne,r}{0}
\node{\cdots}\arrow{e,b}{\pi_{l-2,k}}\arrow{ne,r}{0}
\node{v_{l-2}}\arrow{e,b}{\pi_{l-1,k}}\arrow{ne,r}{0}
\node{v'_{l-1}}\arrow{e,b}{\pi_{l,k}-1}\arrow{ne,r}{2}
\node{v_{l},}
\end{diagram}\nonumber\\
&&(C')={b(v'_l-v_l,\pi_{k,j})\gamma(v_{l-1}-v'_{l-1},\pi_{k,j})
\over \beta(v_{l-1}-v'_{l-1},\pi_{k,j})}\frac{[\pi_{k,j}]}{[\pi_{k,j}+1]}
\times\nonumber\\
&&\begin{diagram}
\node{{v'_{k-1}}}\arrow{e,t}{\pi_{k,j}}
\node{v'_{k}}\arrow{e,t}{\pi_{k+1,j}}
\node{v'_{k+1}}\arrow{e,t}{\pi_{k+2,j}}
\node{\phantom{f}\cdots}\arrow{e,t}{\pi_{l-2,j}}
\node{v'_{l-2}}\arrow{e,t}{\pi_{l-1,j}}
\node{v'_{l-1}}\arrow{e,t}{\pi_{l,k}-1}\arrow{se,r}{2}
\node{v_l}\\
\node[2]{\phantom{f}v_{k}\phantom{f}}\arrow{e,b}{\pi_{k+1,k}}\arrow{ne,r}{0}
\node{v_{k+1}}\arrow{e,b}{\pi_{k+2,k}}\arrow{ne,r}{0}
\node{\cdots}\arrow{e,b}{\pi_{l-2,k}}\arrow{ne,r}{0}
\node{v_{l-2}}\arrow{e,b}{\pi_{l-1,k}}\arrow{ne,r}{0}
\node{v_{l-1}}\arrow{e,b}{\pi_{l,j}-1}
\node{v'_{l}.}
\end{diagram}\nonumber
\ena
Noting that ${h(v_{l-1}-v'_{l-1})}{\beta(v'_{l-1}-v_{l-1},w)}=\beta(v_{l-1}-v'_{l-1},w)$, we can exchange $v_{l-1}$ and $v'_{l-1}$ in $(A')$.
Let $(A'')$ be the term we thus obtain.
Note that $(A')\sim(A'')$. Using the equality
\be
&&f(v'_{l-1}-v_l,2)-\frac{b(v'_l-v_l,w)}{\beta(v_{l-1}-v'_{l-1},w)}\frac{[w]}{[w+1]}f(v_{l-1}-v'_l,2)\nonumber\\
&&={[1][v_l-v'_l+v_{l-1}-v'_{l-1}][v_l-v_{l-1}+{1\over2}][v'_l-v'_{l-1}+{3\over2}]
\over
[v_{l-1}-v'_{l-1}][v'_l-v_l+1][v'_{l-1}-v_l-{1\over2}][v_{l-1}-v'_l-{1\over2}]}
\en
we have
\bea
&&(A)+(A'')={[1][v_l-v'_l+v_{l-1}-v'_{l-1}][v_l-v_{l-1}+{1\over2}][v'_l-v'_{l-1}+{3\over2}]
\over
[v_{l-1}-v'_{l-1}][v'_l-v_l+1][v'_{l-1}-v_l-{1\over2}][v_{l-1}-v'_l-{1\over2}]}\times
\nonumber\\
&&\begin{diagram}
\node{v'_{k-1}}\arrow{e,t}{\pi_{k,j}}
\node{v'_{k}}\arrow{e,t}{\pi_{k+1,j}}
\node{\cdots}\arrow{e,t}{\pi_{l-2,j}}
\node{v'_{l-2}}\arrow{e,t}{\pi_{l-1,j}}
\node{v'_{l-1}}\arrow{e,t}{\pi_{l,j}-1}
\node{v'_l}\\
\node[2]{v_{k}}\arrow{e,b}{\pi_{k+1,k}}\arrow{ne,r}{0}
\node{\cdots}\arrow{e,b}{\pi_{l-2,k}}\arrow{ne,r}{0}
\node{v_{l-2}}\arrow{e,b}{\pi_{l-1,k}}\arrow{ne,r}{0}
\node{v_{l-1}}\arrow{e,b}{\pi_{l,k}-1}
\node{v_l.}
\end{diagram}\lb{b11}
\ena
On the other hand, to calculate $(C)+(C')$, we use the equality
\bea
&&{c(v_{l-1}-v'_{l-1},w_1)b(v'_l-v_l,w_1)\over \beta(v_{l-1}-v'_{l-1},w_1)}
\frac{[w_1]}{[w_1+1]}f(v'_{l-1}-v_l,w_2-1)f(v_{l-1}-v'_l,w_1+w_2-1)
\nonumber\\
&&-{c(v'_l-v_l,w_1)}f(v_{l-1}-v'_l,w_2-1)f(v'_{l-1}-v_l,w_1+w_2-1)
\nonumber\\
&&=-{[1][v_l-v'_l+v_{l-1}-v'_{l-1}][v_{l-1}-v_l+{3\over2}-w_2]
[v'_{l-1}-v'_l+{3\over2}-w_1-w_2]
\over
[v_{l-1}-v'_{l-1}][v'_l-v_l+1][v'_{l-1}-v_l-{1\over2}][v_{l-1}-v'_l-{1\over2}]}.
\nonumber
\ena
Then we have
\bea
&&(C)+(C')
=-{[1][v_l-v'_l+v_{l-1}-v'_{l-1}][v_{l-1}-v_l+{3\over2}-w_2]
[v'_{l-1}-v'_l+{3\over2}-w_1-w_2]
\over
[v_{l-1}-v'_{l-1}][v'_l-v_l+1][v'_{l-1}-v_l-{1\over2}][v_{l-1}-v'_l-{1\over2}]}
\times
\nonumber\\
&&\begin{diagram}
\node{v'_{k-1}}\arrow{e,t}{\pi_{k,j}}
\node{v'_{k}}\arrow{e,t}{\pi_{k+1,j}}
\node{\cdots}\arrow{e,t}{\pi_{l-2,j}}
\node{v'_{l-2}}\arrow{e,t}{\pi_{l-1,j}}
\node{v'_{l-1}}\arrow{e,t}{2}
\node{v'_l}\\
\node[2]{v_{k}}\arrow{e,b}{\pi_{k+1,k}}\arrow{ne,r}{0}
\node{\cdots}\arrow{e,b}{\pi_{l-2,k}}\arrow{ne,r}{0}
\node{v_{l-2}}\arrow{e,b}{\pi_{l-1,k}}\arrow{ne,r}{0}
\node{v_{l-1}.}
\end{diagram}\lb{b12}
\ena
Comparing \eqref{b11} and \eqref{b12}, we have $(A)+(A'')+(C)+(C')=0$.


\section{$RLL=LLR^*$ Relation}\lb{RLL}
We here derive some of the relations of the 
half currents involved in the $RLL$-relation \eqref{thm:RLL},
and compare them with those in Theorem \ref{KEF}. 

From the definition \eqref{def:lhat}, the components
of the $L$-operator $\widehat{L}^+(v)$ are given by 
\begin{eqnarray}
&&\widehat{L}^{+}_{ll}(v)=K^+_l(v)+\sum_{m=1}^N F^+_{l,m}(v)K^+_m(v)
E^+_{m,l},\\
&&\widehat{L}^+_{k l}(v)=\left\{\begin{array}{cc}
F_{k,l}^+(v)K_l^+(v)+\sum_{m=l+1}^N F^+_{k,m}(v)K^+_m(v)E^+_{m,l}& (k<l),\\
K_k^+(v)E_{k,l}^+(v)+\sum_{m=k+1}^N F^+_{k,m}(v)K^+_m(v)E^+_{m,l}& (k>l).
\end{array}\right.
\end{eqnarray}
It is convenient to introduce the reduced $R$-matrix 
and $L$-operators,  $R^+(v,s|j)$ and 
$\widehat{L}^+(v|j) \ (1\leq j\leq N)$,  by
\begin{eqnarray}
&&R^+(v,s|j)=\left(\ R^{+ m n}_{\ \ k l}(v,s)\ 
\right)_{j\leq k,l,m,n \leq N},\\
&&\widehat{L}^+(v|j)=\left(\ \widehat{L}^+_{k l}(v)\ \right)_{j\leq k,l \leq N}.\lb{reducedl}
\end{eqnarray}
Then the inverse of $L^+(v|j)$ is given by
\begin{eqnarray}
&&L^+(v|j)^{-1}\nonumber\\
&=&\left(
\begin{array}{cccc}
K_j^{+-1}&-K_j^{+-1}F_{j,j+1}^+&K_j^{-1}x_j&*\\
-E_{j+1,j}^+K_j^{+-1}&
E_{j+1,j}^+K_j^{+ -1}F_{j,j+1}^++K_{j+1}^{+ -1}
&
-E_{j+1,j}^+K_j^{+ -1}x_j-K_{j+1}^{+ -1}F_{j+1,j+2}^+
&*\\
y_jK_j^{-1}&-y_jK_j^{+ -1}F_{j,j+1}^+-E_{j+2,j+1}^+ K_{j+1}^{+ -1}&*&* \\
\ast &* &* &* \\
\end{array}
\right).\nonumber\\
\end{eqnarray}
Here we omitted the argument $v$ and set
\begin{eqnarray}
x_j(v)=F_{j,j+1}^+(v)F_{j+1,j+2}^+(v)-F_{j,j+2}^+(v),\\
y_j(v)=E_{j+2,j+1}^+(v)E_{j+1,j}^+(v)-E_{j+2,j}^+(v).
\end{eqnarray}

Due to the speciality of the form of the $R$-matrix \eqref{rmat}, 
we have the reduced relation
\bea
&&R^{+(1,2)}(v,P+h|j)L^{+(1)}(v_1|j)L^{+(2)}(v_2|j)=L^{+(2)}(v_2|j)L^{+(1)}(v_1|j)R^{*+}(v,P|j).
\lb{redRLL}
\ena
In the below, we  use this rather in its inverted form   
\bea
&&L^{+(1)}(v_1|j)^{-1}L^{+(2)}(v_2|j)^{-1}R^{+(1,2)}(v,P+h|j)=R^{*+(1,2)}(v,P|j)L^{+(2)}(v_2|j)^{-1}
L^{+(1)}(v_1|j)^{-1},\nn\\
\lb{redinvRLL1}\\
&&L^{+(2)}(v_2|j)^{-1}R^{+(1,2)}(v,P+h|j)L^{+(1)}(v_1|j)=L^{+(1)}(v_1|j)R^{*+(1,2)}(v,P|j)L^{+(2)}(v_2|j)^{-1}
.\nn\\
\lb{redinvRLL2}
\ena

\subsection{Relations among $K^+_j(v)$'s}
Now some of the relations among {$K^+_j(v)\ (1\leq j\leq N)$ } are derived as
 follows.
The $(N,N),(N,N)$ component of the $RLL=LLR^*$ relation \eqref{thm:RLL} 
yields 
\begin{eqnarray}
K_N^+(v_1)K_N^+(v_2)
=\rho(v_1-v_2)K_N^+(v_2)K_N^+(v_1). 
\end{eqnarray}
Similarly, the $(j,j),(j,j)$ component of the $L_j^{-1}L_j^{-1}R=R^*L_j^{-1}L_j^{-1}$ relation (\ref{redinvRLL1})\ $(1\leq j \leq N-1)$ yields
\begin{eqnarray}
K_j^+(v_1)K_j^+(v_2)
=\rho(v_1-v_2)K_j^+(v_2)K_j^+(v_1) 
\end{eqnarray}
and the $(N,j),(N,j)$ component of the $L_j^{-1} R L_j=L_j R^*L_j^{-1}$
 relation  (\ref{redinvRLL1})
\ $(1\leq j \leq N-1)$ yields
\begin{eqnarray}
K_j^+(v_1)K_N^+(v_2)
=\rho(v_1-v_2)\frac{[v_1-v_2-1]^*[v_1-v_2]}{[v_1-v_2]^*[v_1-v_2-1]}
K_N^+(v_2)K_j^+(v_1).
\end{eqnarray}
These relations coincide with the relations \eqref{KjKj} and
\eqref{kjkl}.

\subsection{Relations between $K^+_N(v)$ and $E^+_{N,j}(v)$}
The  $(N,N),(N,j)$  components of the $ RLL=LLR^*$ relation 
\eqref{thm:RLL}\ $(1\leq j \leq N-1)$ yields
\begin{eqnarray}
K_N^+(v_1)^{-1}E_{N,j}^+(v_2)
K_N^+(v_1)=E_{N,j}^+(v_2)\frac{1}{\bar{b}_{}^*(v_1-v_2)}-
E_{n,j}^+(v_1)\frac{c_{}^*(v_1-v_2,P_{j,N})}{\bar{b}_{}^*(v_1-v_2)}.
\end{eqnarray}
This coincides with the case $l=N$ of \eqref{klelj}.


\subsection{Relations between $K^+_N(v)$ and $F^+_{j,N}(v)$}
The $(N,j),(N,N)$  components of the $RLL=LLR^*$ relation 
\eqref{thm:RLL}\ $(1\leq j \leq N-1)$ yields
\begin{eqnarray}
K_N^+(v_1)F_{j,N}^+(v_2)
K_N^+(v_1)^{-1}=
\frac{1}{\bar{b}_{}(v_1-v_2)}F_{j,N}^+(v_2)-
\frac{\bar{c}_{}(v_1-v_2,P_{j,N}+h_{j,N})}{
\bar{b}_{}(v_1-v_2)}F_{j,N}^+(v_1).
\end{eqnarray}
This coincides with the case $l=N$ of \eqref{klfjl}.


\subsection{Relations among $E^+_{l,j}(v)$'s}
The $(N,N),(j,j)$  component of the $RLL=LLR^*$ relation 
\eqref{thm:RLL}\ $(1\leq j \leq N-1)$ yields 
\begin{eqnarray}
K_N^+(v_1)E_{N,j}^+(v_1)
K_N^+(v_2)E_{N,j}^+(v_2)=
\rho(v_1-v_2)K_N^+(v_2)E_{N,j}^+(v_2)
K_N^+(v_1)E_{N,j}^+(v_1).\lb{ee1}
\end{eqnarray}
The $(N,N),(k,j)$  component of the $RLL=LLR^*$ relation 
\eqref{thm:RLL}\ $(1\leq j,k \leq N-1,\ j\not=k )$ yields 
\bea
&&\rho^+(v_1-v_2)K_N^+(v_1)E^+_{N,k}(v_1)K_N^+(v_2)E^+_{N,j}(v_2)\nn\\
&&\qquad\qquad=K_N^+(v_2)E^+_{N,j}(v_2)K_N^+(v_1)E^+_{N,j}(v_1)
{R}^{*kj}_{kj}(v_1-v_2,P_{j,k})\nn\\
&&\qquad\qquad
+K_N^+(v_2)E^+_{N,j}(v_2)K_N^+(v_1)E^+_{N,k}(v_1)
{R}^{*kj}_{jk}(v_1-v_2,P_{j,k}).\lb{ee2}
\ena
After a little calculation using \eqref{KjKj}, \eqref{ee1} and \eqref{ee2} coincide with the case $l=N$ in 
\eqref{eljelj2} and \eqref{elkelj}, respectively.

\subsection{Relations among $F^+_{j,l}(v)$'s}
The $(j,j),(N,N)$ component
 of the $RLL=LLR^*$ relation \eqref{thm:RLL}\ $(1\leq j \leq N-1)$ yields
\begin{eqnarray}
F_{j,N}^+(v_1)K_N^+(v_1)
F_{j,N}^+(v_2)K_N^+(v_2)=
\rho(v_1-v_2)
F_{j,N}^+(v_2)K_N^+(v_2)
F_{j,N}^+(v_1)K_N^+(v_1).\lb{ff1}
\end{eqnarray}
The $(j,k), (N,N)$  component of the $RLL=LLR^*$ relation 
\eqref{thm:RLL}\ $(1\leq j,k \leq N-1,\ j\not=k )$ yields 
\bea
&&\rho^{+*}(v_1-v_2)F^+_{k,N}(v_2)K_N^+(v_2)F^+_{j,N}(v_1)K_N^+(v_1)\nn\\
&&\qquad\qquad={R}^{jk}_{jk}(v,P_{j,k}+h_{j,k})F^+_{j,N}(v_1)K_N^+(v_1)
F^+_{k,N}(v_2)K_N^+(v_2)\nn\\
&&\qquad\qquad
+{R}^{kj}_{jk}(v,P_{j,k}+h_{j,k})F^+_{k,N}(v_1)K_N^+(v_1)F^+_{j,N}(v_2)
K_N^+(v_2).\lb{ff2}
\ena
After a little calculation using \eqref{KjKj}, 
\eqref{ff1} and \eqref{ff2} coincide with the case $l=N$ in 
\eqref{fjlfjl2} and \eqref{fklfjl}, respectively.

\subsection{The relations between $E_{l,j}^+$'s and $F^+_{k,l}$'s}
The 
$(j,N),(N,N-1)$ component together with the $(j,N),(N,N)$ and 
$(N,N),(N,N-1)$ components of the $RLL=LLR^*$ relation \eqref{thm:RLL} 
$(1\leq j \leq N-2)$ yield
\begin{eqnarray}
&&[E_{N,N-1}^+(v_2),F_{j,N}^+(v_1)]=
K_N^+(v_2)^{-1}\frac{c(v_1-v_2,P_{j,N}+h_{j,N})}{
\bar{b}(v_1-v_2)}F_{j,N-1}^+(v_2)K_{N-1}^+(v_2)\nn\\
&&\qquad\qquad\qquad\qquad\qquad
-F_{j,N-1}^+(v_1)K_{N-1}^+(v_1)\frac{c^*(v_1-v_2,P_{N-1,N})}
{\bar{b}^*(v_1-v_2)}K_N^+(v_1)^{-1}.\lb{ef1}
\end{eqnarray}
The $(j+1,j),(j,j+1)$ component 
together with the 
$(j+1,j),(j,j)$ and
$(j,j),(j,j+1)$ components of the 
$ L_j^{-1}L_j^{-1}R=R^* L_j^{-1}L_j^{-1}$ relation \eqref{redinvRLL1}
$ (1\leq j \leq N-1)$ yield
\begin{eqnarray}
&&[E_{j+1,j}^+(v_1),F_{j,j+1}^+(v_2)]
=K_j^+(v_2)\frac{\bar{c}^*(v_1-v_2,P_{j,j+1})}{
\bar{b}^*(v_1-v_2)}K_{j+1}^+(v_2)^{-1}\nn\\
&&\qquad\qquad\qquad\qquad\qquad
-K_{j+1}^+(v_1)^{-1}
\frac{\bar{c}(v_1-v_2,P_{j,j+1}+h_{j,j+1})}{\bar{b}(v_1-v_2)}
K_j^+(v_1).\lb{ef2}
\end{eqnarray}
These equations \eqref{ef1} and \eqref{ef2} coincide with the cases $l=N$ and 
$l=j+1$ of \eqref{elfjl}, respectively. 

Similarly, the $(N-1,N),(N,j)$ component together with the 
 $(N-1,N),(N,N)$ and $(N,N),(N,j)$ components of the 
 $RLL=LLR^*$ relation \eqref{thm:RLL} $(1\leq j \leq N-2)$ yield
\begin{eqnarray}
&&[E_{N,j}^+(v_2),F_{N-1,N}^+(v_1)]
=K_N^+(v_2)^{-1}\frac{c(v_1-v_2,P_{N-1,N}+h_{N-1,N})}
{\bar{b}(v_1-v_2)}K_{N-1}^+(v_2)E_{N-1,j}^+(v_2)\nn\\
&&\qquad\qquad\qquad\qquad\qquad-
K_{N-1}^+(v_1)E_{N-1,j}^+(v_1)
\frac{c^*(v_1-v_2,P_{j,N})}{\bar{b}^*(v_1-v_2)}K_N^+(v_1)^{-1}.\lb{ef3}
\end{eqnarray}
The $(j,j+1),(j+1,j)$
component together with the  $(j,j),(j+1,j)$
and $(j,j+1),(j,j)$ components of the $ L_j^{-1}L_j^{-1}R=R^* L_j^{-1}L_j^{-1}$
 relation \eqref{redinvRLL1}$ (1\leq j \leq N-1)$ yield
\begin{eqnarray}
&&[E_{j+1,j}^+(v_2),F_{j,j+1}^+(v_1)]=
K_{j+1}^+(v_2)^{-1}\frac{c(v_1-v_2,P_{j,j+1}+h_{j,j+1})}{
\bar{b}(v_1-v_2)}K_j^+(v_2)\nn\\
&&\qquad\qquad\qquad\qquad\qquad
-K_{j}^+(v_1)\frac{c^*(v_1-v_2,P_{j,j+1})}{
\bar{b}^*(v_1-v_2)}K_{j+1}^+(v_1)^{-1}.\lb{ef4}
\end{eqnarray}
These equations \eqref{ef3} and \eqref{ef4} coincide with the cases $l=N$ and 
$l=j+1$ of \eqref{eljfl}, respectively. 

Finally, the following relations with $j\leq N-2$ are examples
of those which we have not yet checked for our half currents.
\begin{eqnarray}
&&[E_{N,j}^+(v_2),F_{j,N}^+(v_1)]\nn\\
&&=
K_N^+(v_2)^{-1}\frac{c(v_1-v_2,P_{j,N}+h_{j,N})}{
\bar{b}(v_1-v_2)}K_j^+(v_2)-
K_j^+(v_1)\frac{c^*(v_1-v_2,P_{j,N})}{
\bar{b}^*(v_1-v_2)}K_N^+(v_1)^{-1}\nonumber\\
&&+
\sum_{k=j+1}^{N-1}
\left(K_N^+(v_2)^{-1}
\frac{c(v_1-v_2,P_{j,N}+h_{j,N})}{\bar{b}(v_1-v_2)}
F_{j,k}^+(v_2)K_k^+(v_2)E_{k,j}^+(v_2) \right. \nn\\
&&\qquad\qquad\qquad\qquad \left. -
F_{j,k}^+(v_1)K_k^+(v_1)E_{k,j}^+(v_1)
\frac{c^*(v_1-v_2,P_{j,N})}{
\bar{b}^*(v_1-v_2)}K_N^+(v_1)^{-1}
\right).
\end{eqnarray}
These are derived from the $(j,N),(N,j)$ components together 
with the $(j,N),(N,N)$ 
and $(N,N),(N,j)$ components of the $RLL=LLR^*$ relation 
\eqref{thm:RLL} $(1\leq j \leq N-1)$. 

\section{Evaluation Module}\lb{Evaluation}

We here summarize the evaluation module $(\pi_{V,z},V_z=V[z,z^{-1}])$ of 
$U_q(\widehat{\goth{sl}}_N)$ 
associated with the vector representation $
V={\mathbb{C}}^N$.

The evaluation module $(\pi_z, V_z)$ in terms of the Drinfeld generators,
is defined by the following formulae.
\begin{eqnarray}
&&\pi_z(c)=0,~\pi_z(d)=z\frac{d}{dz},\\
&&\pi_z(a_{j,n})=\frac{[n]}{n}(q^{j-N+1}z)^n
(q^{-n}E_{j j}-q^n E_{j+1 j+1}),\\
&&\pi_z(x_{j,n}^+)=(q^{j-N+1} z)^n E_{j j+1},\\
&&\pi_z(x_{j,n}^-)=(q^{j-N+1} z)^n E_{j+1 j},\\
&&\pi_z(h_j)=E_{j j}-E_{j+1 j+1},~ 
\pi_z(h_{\bar{\epsilon}_j})=-E_{j j}.
\end{eqnarray}
Then the elliptic currents 
$k_j(w,p), \psi_j^\pm(w,p), e_j(w,p), f_j(w,p)$ of $U_{q}(\slnh)$ defined in
\eqref{def:e}-\eqref{def:psi-} are represented by
\begin{eqnarray}
&&\pi_{z}(k_j(w,p))=
\frac{\{q^{r+2N+1}z/w\}
\{q^{r+1}z/w\}
\{q^{r-1}w/z\}
\{q^{r-2N+3}w/z\}
}{
\{q^{r+2N-1}z/w\}
\{q^{r+3}z/w\}
\{q^{r+1}w/z\}
\{q^{r-2N+1}w/z\}
}\nonumber\\
&&\qquad\qquad\qquad\times
\left(\frac{\Theta_p(q^{r+1}z/w)}{\Theta_p(q^{r-1}z/w)}
\sum_{k=1}^{j-1}E_{k k}+
E_{j j}+
\frac{\Theta_p(q^{r+3}z/w)}{\Theta_p(q^{r+1}z/w)}
\sum_{k=j+1}^N E_{k k}\right),\nonumber\\
\\
&&\pi_{z}(\psi_j^\pm(q^{\mp r}w,p))=q^{\pm h_j}
\frac{\Theta_p(q^{r-j+2h_j+N-1}w/z)}
{\Theta_p(q^{r-j+N-1}w/z)},\\
&&\pi_{z}(e_j(w,p))=E_{j j+1} 
\frac{(pq^{2};p)_\infty}
{(p;p)_\infty} \delta(q^{j-N+1}z/w),\\
&&\pi_{z}(f_j(w,p))=E_{j+1 j} \frac{(pq^{-2};p)_\infty}
{(p;p)_\infty} 
\delta(q^{j-N+1}z/w),
\end{eqnarray}
where $\{z\}=(z;p,q^{2N})_\infty$.
Especially, the auxiliary currents $u_j^\pm(w,p)$ are represented by
\begin{eqnarray}
\pi_z(u_j^+(w,p))=\frac{(pq^{-j+2h_j+N-1}w/z;p)_\infty}{
(pq^{-j+N-1}w/z;p)_\infty},~~
\pi_z(u_j^-(w,p))=
\frac{(pq^{j-2h_j-N+1}z/w;p)_\infty}{
(pq^{j-N+1}z/w;p)_\infty}.
\end{eqnarray}
Due to this representation, we can obtain
the representation  of the half currents.
After getting rid of some unpleasant fractional power factors of
$q$ and $z$ by a certain gauge transformation, we have
the following result.
\begin{eqnarray}
&&\pi_{v_2}(K_j^+(v_1)e^{-Q_{\bar{\epsilon}_j}})
=\rho^+\left(v_1-v_2\right)
\nonumber\\
&&\qquad\qquad\qquad \times
\left(\frac{[v_1-v_2]}{
[v_1-v_2+1]}
\sum_{k=1}^{j-1}E_{k k}+
E_{j j}+
\frac{[v_1-v_2-1]}{[v_1-v_2]}
\sum_{k=j+1}^N E_{k k}\right),\nn\\
&&\\
&&\pi_{v_2}(e^{-\eta_j}F_{j,l}^+(v_1)e^{\eta_l})=
E_{l j}\frac{[v_1-v_2+P_{j,l}-1][1]}
{[v_1-v_2][P_{j,l}-1]}
,\\
&&\pi_{v_2}(e^{Q_{\bar{\epsilon}_{l}-\eta_l}}
E_{l,j}^+(v_1)e^{-Q_{\bar{\epsilon}_{j}+\eta_j}}
)=-E_{j l}\frac{[v_1-v_2-P_{j,l}][1]}
{[v_1-v_2][P_{j,l}]},
\end{eqnarray}
where $z=q^{2v}$.
It is easy to check that these quantities satisfy
the commutation relations of the half currents
$K_j^+(v), F_{j,l}^+(v)$ and $ E_{l,j}^+(v)$.

Finally, let us check the results by calculating 
the $R$-matrix as the image of
the $L$-operator $\hat{L}^+(v)$ in \eqref{def:lhat}.
\begin{eqnarray}
R^+(v_1-v_2,P)=(\pi_{v_2} \otimes id)\hat{L}^+(v_1).
\end{eqnarray}
Using \eqref{kjpjl} and Riemann's theta identity,  
we obtain the following.
\begin{eqnarray}
R^+(v,P)=
\rho^+\left(v\right)
\left(\begin{array}{cccc}
R_{11}(v,P)&\cdots&\cdots&R_{1N}(v,P)\\
R_{21}(v,P)&\cdots&\cdots&R_{2N}(v,P)\\
\vdots&&\ddots&\vdots\\
R_{N1}(v,P)&\cdots&\cdots&R_{NN}(v,P)
\end{array}\right),
\end{eqnarray}
where
\begin{eqnarray}
&&R_{j j}(v,P)=
\sum_{k=1}^{j-1}
\bar{b}\left(v\right)
E_{k k}+E_{j j}+
\sum_{k=j+1}^{N}
b\left(v,P_{j,k}
\right)E_{k k}~~\qquad(1\leq j \leq N),\nonumber\\
\\
&&R_{jl}(v,P)=
c\left(v,P_{j,l}\right)E_{l j}~~
,\\
&&
R_{l j}(v,P)=
\bar{c}\left(v,P_{j,l}\right)E_{j l}~~\qquad(1\leq j<l \leq N).
\end{eqnarray}
This expression coincides with the $R$-matrix 
given by (\ref{rmatfull}).

\end{appendix}


\end{document}